\documentclass[a4paper,oneside,11pt]{article}
\usepackage{amsmath}
\usepackage{amsmath}
\usepackage{amsfonts}
\usepackage{amssymb}
\usepackage{titlesec}
\usepackage{graphicx}
\usepackage{subfigure}
\usepackage{soul}
\usepackage{fancyhdr}
\usepackage{float}
\usepackage{pdfsync}
\usepackage{latexsym}
\usepackage{mathrsfs}
\usepackage{mathrsfs}
\usepackage{mathtools}
\usepackage{amsthm}
\usepackage{wasysym}
\usepackage{cite}
\usepackage{color}
\usepackage{geometry}
\usepackage{extarrows}
\usepackage{enumerate}
\usepackage{dsfont}
\usepackage{bm}
\usepackage[marginal]{footmisc}
\usepackage[dvips]{epsfig}
\usepackage{amscd}
\usepackage[marginal]{footmisc}
\usepackage[bookmarks,colorlinks,citecolor=blue,linkcolor=red]{hyperref}
\allowdisplaybreaks

\numberwithin{equation}{section}

\newtheorem{theorem}{Theorem}[section]
\newtheorem{lemma}{Lemma}[section]
\newtheorem{remark}{Remark}[section]
\newtheorem{proposition}{Proposition}[section]
\renewcommand{\thefootnote}{}

\title{Phase Transition in Non-isentropic Compressible  Immiscible Two-Phase Flow with van der Waals Equation of State  }
\date{  }
\author{Yazhou C{\small HEN}$^1$, Yi P{\small ENG}$^1$, Xiaoding S{\small HI}$^{1*}$, Xiaoping W{\small ANG}$^{2,3}$ \\[3mm]
\scriptsize$^{1}$ {College of Mathematics and Physics, Beijing University of
Chemical Technology, Beijing 100029, China}\\
\scriptsize$^2$ {School of Science and Engineering, The Chinese University of Hong Kong, Shenzhen, Guangdong 518172, China}\\
\scriptsize$^3$ {Shenzhen International Center for Industrial and Applied Mathematics, Shenzhen Research Institute of Big Data,}\\
\scriptsize{ Guangdong 518172,China  }
}

\begin{document}
\maketitle
\renewcommand{\thefootnote}{\fnsymbol{footnote}}
\footnotetext[1]{{Corresponding author. }\\
{Email: chenyz@mail.buct.edu.cn (Y.Chen), 2024500077@buct.edu.cn (Y.Peng), shixd@mail.buct.edu.cn (X.Shi), wangxiaoping@cuhk.edu.cn, (X.Wang)}}

\begin{abstract}
This study establishes the global well-posedness of the compressible non-isentropic Navier-Stokes/Allen-Cahn system governed by the van der Waals equation of state $p(\rho,\theta)=- a\rho^2+\frac{R\theta\rho}{1-b\rho}$ and degenerate thermal conductivity $\kappa(\theta)=\tilde{\kappa}\theta^\beta$, where $p$, $\rho$ and $\theta$  are the pressure, the density and the temperature of the flow respectively, and $a,b,R,\tilde\kappa$ are positive constants related to the physical properties of the flow.  Navier-Stokes/Allen-Cahn system models immiscible two-phase flow with diffusive interfaces, where the non-monotonic pressure-density relationship in the van der Waals equation drives gas-liquid phase transitions. By developing a refined $L^2$-energy framework, we prove the existence and uniqueness of global strong solutions to the one-dimensional Cauchy problem  for non-vacuum and finite-temperature  initial data, without imposing smallness restrictions on the initial conditions. The findings demonstrate that despite non-monotonic pressure inducing substantial density fluctuations and triggering phase transitions, all physical quantities remain bounded over finite time intervals.
\end{abstract} 

\


\section{Introduction and Main Theorems}\label{sec:int}

\ \ \ \ 
Phase transition in immiscible two-phase flow is a complex intersection involving fluid mechanics, thermodynamics and mass and heat transfer, which is commonly found in energy, chemical, aerospace and other fields (such as oil and gas extraction, nuclear reactor cooling, refrigeration systems, etc.). The problem of gas-liquid phase transition is closely related to van der Waals equation of state proposed by  van der Waals \cite{W-1894}. By introducing the modified ideal gas law of intermolecular force and molecular volume, van der Waals equation can more accurately describe the critical behavior of gas-liquid transformation, gas-liquid coexistence zone and non-ideal characteristics of phase transition.
The van der Waals equation was the first successful equation of state to describe the gas-liquid phase transition, and its form is:
\begin{equation}\label{Original van der Waals equation}
 (p+\frac{a}{v^2})(v-b)=R\theta,
\end{equation}
where $p$, $v$ and $\theta$  are the pressure, the volume and the temperature of the flow respectively, $R$ is the gas constant,  $a$ and $b$ are constants characterizing the effect of the molecular cohesive forces and the finite size of the molecules. The internal pressure term $\frac{a}{v^2}$ reflects the correction of pressure due to intermolecular attractions, and directly affects the potential energy component of the internal energy. The volume correction term $b$ takes into account the volume occupied by the molecules themselves, and indirectly influences the calculation of the internal energy. According to the thermodynamic relationship, the volume dependence of internal energy $E_{int}$ can be expressed as $\frac{\partial E_{int}}{\partial v}=\frac{a}{v^2}$, after integrating, the expression for internal energy is obtained
\begin{equation}\label{internal energy}
  E_{int}=c_v\theta-\frac a{v}+E_{int_0}.
\end{equation}
Here, $c_v$ represents the constant volume heat capacity, and $E_{int_0}$ is the integration constant.


Formally, as $v$ approaches $b^+$ (meaning that the molecules are packed so densely that no further compression is possible), according to the mathematical derivation of the van der Waals equation, either $p$ should approach $+\infty$ or $\theta$ should approach $0^+$. However, in reality, neither of these scenarios can occur. The underlying reason is that when the specific volume 
$v$ is compressed to a value close to the molecular diameter, quantum effects inevitably emerge. These quantum effects manifest in various forms, such as the repulsion between electron clouds, the discretization of energy levels, and changes in statistical behavior. Therefore, in this paper, to circumvent the quantum effects that occur near the molecular diameter $b^+$, we propose a modified van der Waals relation:


\

\noindent\textbf{van der Waals Equation of State.}
\begin{equation}\label{modified van der Waals equation}
p(v,\theta)=\left\{\begin{array}{llll}
\displaystyle -\frac a{v^2}+\frac{R\theta}{v-b},&\displaystyle v>b+h,\\
\displaystyle +\infty,&\displaystyle \mathrm{other},
\end{array}\right.
\end{equation}
where, physically speaking, $h>0$ is actually a positive constant small enough and it is related to the specific material properties.

\

\begin{remark}\label{The positive nature of pressure}
From the expression \eqref{modified van der Waals equation}, one knows that the isotherm of the van der Waals equation \eqref{modified van der Waals equation} will appear as a "hump region" on the $p-v$ diagram (see Figure 1(a))  as $0<\theta<\frac{8a}{27Rb}$, i.e. there is an interval of $v$ where $p$ increases as $v$ increases, and the phase transition will take place around this region. Moreover, \eqref{modified van der Waals equation} ensures the positivity of the pressure. In fact, under the condition \eqref{modified van der Waals equation}, by adding a sufficiently large positive constant, the pressure can be made to take on a positive value.Therefore, without loss of generality, throughout this paper, we always assume the positive condition of pressure, namely
\begin{equation}\label{Hypothesis of p}
  p(v,\theta)>0.
\end{equation}
\end{remark}

\

\begin{lemma}\label{properities of p}
 The pressure $p$ has the following properties (see Figure 1(b))
  \begin{enumerate}
  \item[1]  $\displaystyle\lim_{v\rightarrow+\infty}p(v,\theta)=0$;

\item[2] For a given isotherm (fixed temperature $\theta\in(0,\theta_c)$, $\theta_c$ is the critical temperature of the fluid, see figure 1(a)), the pressure $p$ exhibits non-monotonic dependence on volume $v$, characterized by the existence of two critical points: a local minimum at $v_\alpha(\theta)$ and a local maximum at $v_\beta(\theta)$, and $b+h<v_\alpha(\theta)<v_\beta(\theta)$, (with $v_\alpha(\theta) < v_\beta(\theta)$). The function $p$ displays the following monotonic behavior:
\begin{enumerate}    
\item[(i)] Decreasing domains (called as hyperbolic domain): the pressure $p$ decreases monotonically with respective of $v$ in the intervals ($b+h,v_\alpha(\theta)$) and ($v_\beta(\theta),+\infty$) for fixed $\theta$.

\item[(ii)] Increasing domain (called as elliptic domain): the pressure  $p$ increases monotonically with respective of $v$ in the intervals ($v_\alpha(\theta),v_\beta(\theta)$)  for fixed $\theta$.
 
\item[(iii)] Additionally, there are two unique points $v_*(\theta)$ and $v^*(\theta)$, satisfying $v_* (\theta)< v_\alpha(\theta)$, and $v^*(\theta) > v_\beta(\theta)$. These points define a horizontal line segment between $B$ and $C$ on the $p$-isotherm in Fig.1(b). This line corresponds to a specific pressure value where the areas labeled I ($S_{CGF}$) and II ( $S_{GEB}$) in Fig. 1(b) are equal, a condition central to the \textbf{Maxwell equal area rule}(see Remark \ref{Maxwell}).
  \end{enumerate} 

 \item[3] The partial derivatives of pressure with respect to volume $v$ and temperature $\theta$ are given by:
 \begin{equation}\label{the derivative of pressure v}
  \frac{\partial p}{\partial v}=\left\{\begin{array}{llll}
  \displaystyle\frac {2a}{v^3}-\frac{R\theta}{(v-b)^2},&v>b+h,\\
 \displaystyle +\infty,&\mathrm{ other},
  \end{array}\right.\ \frac{\partial p}{\partial \theta}=\left\{\begin{array}{llll}
  \displaystyle\frac{R}{v-b},&v>b+h,\\
 \displaystyle +\infty,& \mathrm{other}.
  \end{array}\right.
  \end{equation}
\end{enumerate}
\end{lemma}

 \begin{figure}
\centering
\subfigure[Figure of $p(v,\theta)$ isotherm at various $\theta$]{
\begin{minipage}[t]{0.45\textwidth}
\centering
\includegraphics[width=2.8 in,height=2.2 in]{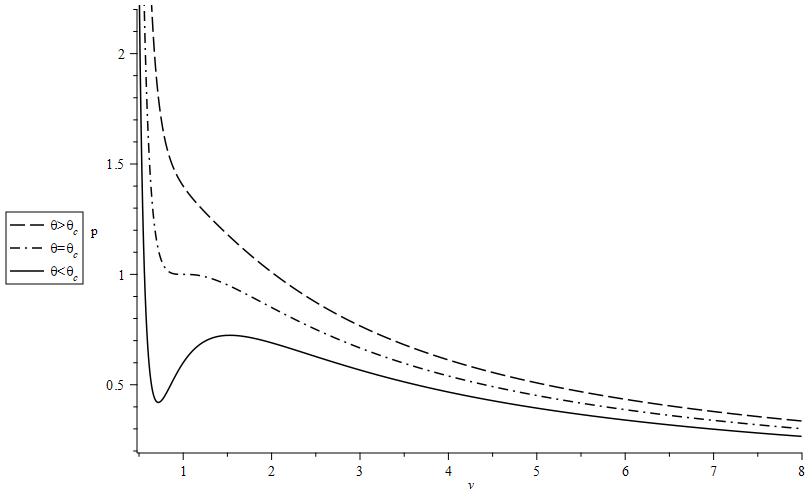}
\end{minipage}}
\hfill
\subfigure[$p-v$ isotherm  with subcritical $\theta$, ($0<\theta<\frac{8a}{27Rb}$) ]{
\begin{minipage}[t]{0.45\textwidth}
\centering
\includegraphics[width=2.8 in,height=2.2 in]{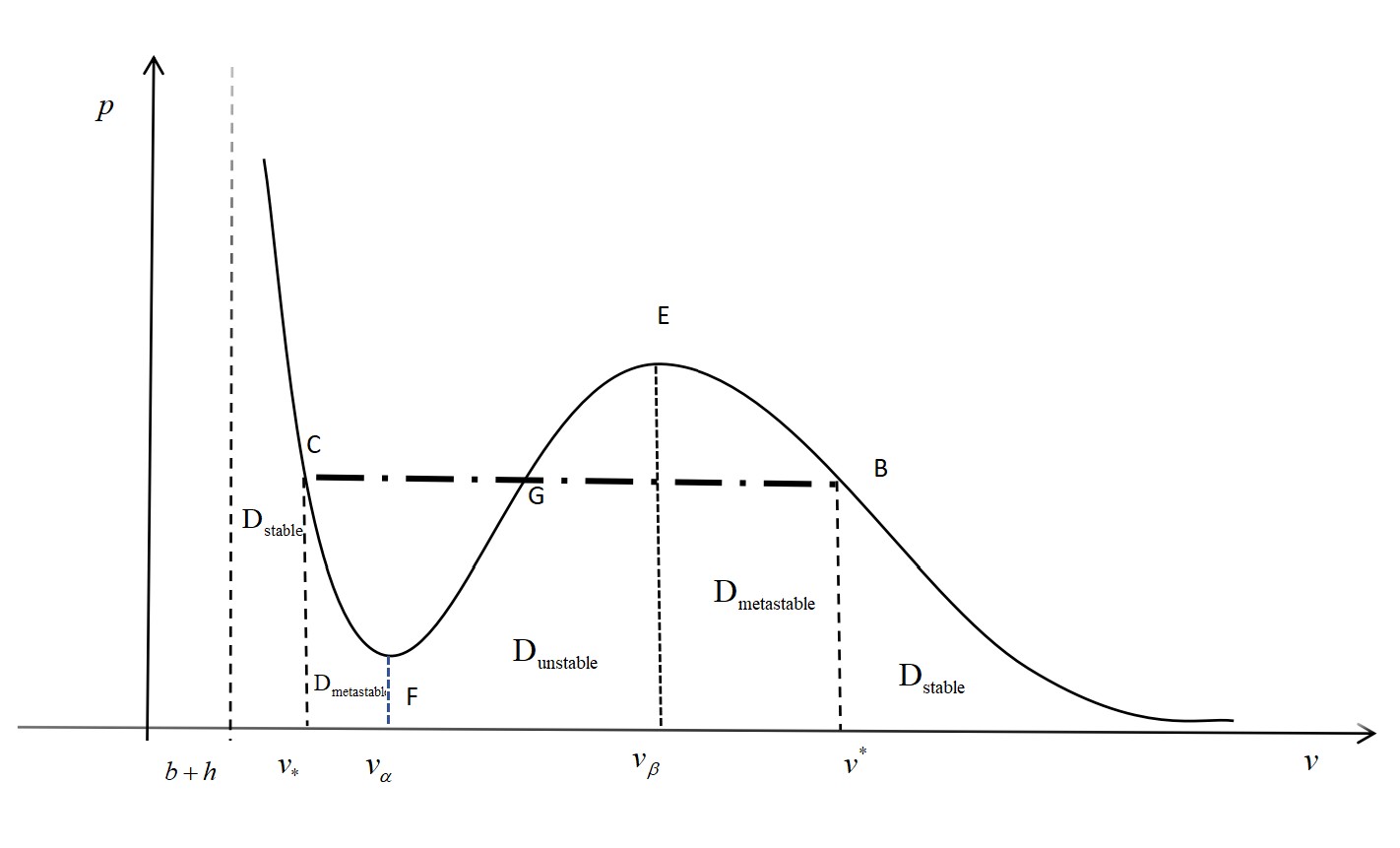}
\end{minipage}}
\centering
\caption{Figure of $p-v$}
\label{model2figure1}
\end{figure}
\begin{remark}\label{Maxwell}
\textbf{(Maxwell equal area rule.)}
In the theoretical framework of thermodynamic phase transitions, the  Maxwell equal area rule (or isobaric equilibrium lines) stands as a cornerstone for elucidating coexistence phenomena, such as gas-liquid and solid-liquid phase transitions. Its fundamental principle lies in enforcing pressure equilibrium across coexisting phases by mathematically imposing that the areas under the Gibbs free energy curves for the two phases must be equal at a given temperature, i.e.,  the areas $S_{CGF}$ and $S_{GEB}$  under the isotherm must satisfy the condition of thermodynamic equilibrium, requiring the areas above and below the horizontal segment to be equal. This geometric constraint ensures thermodynamic consistency. 
\end{remark}

\

After the above historical review and analysis, we know that the van der Waals equation refines the ideal gas law by incorporating intermolecular forces and accounting for the finite volume of molecules. Consequently, it is capable of describing phase transitions and critical phenomena in the vicinity of the liquid-gas coexistence curve, thereby establishing itself as a foundational model in thermodynamics and statistical mechanics.
In this paper, based on the van der Waals equation of state, we investigate the phase transition problem in immiscible two-phase flow.
An important and widely adopted model for characterizing immiscible two-phase flow is the Navier-Stokes/Allen-Cahn system. We shall study the well-posedness  for compressible non-isothermal Navier-Stokes/Allen-Cahn system  with van der Waals equation of state. We are particularly interested in development of the solution with initial data inside the elliptic region. The Navier-Stokes/Allen-Cahn system is one of the important mathematical models for describing immiscible two-phase flow with diffusion interface (Blesgen \cite{B1999}, Heida-M$\mathrm{\acute{a}}$lek-Rajagopal  \cite{HMR2012},  for details), which is given by
\begin{equation}\label{NSFAC3d}
\left\{\begin{array}{llll}
\displaystyle \rho_t+\textrm{div}(\rho \mathbf{u})=0,\\
\displaystyle (\rho \mathbf{u})_{t}+\mathrm{div}\big(\rho \mathbf{u}\otimes\mathbf{u}\big)-\nu\Delta\mathbf{u}-\big(\nu+\lambda\big)\nabla\mathrm{div}\mathbf{u}\\
\displaystyle\qquad\qquad\qquad\qquad\qquad=\mathrm{div}\big(\frac{\epsilon}{2}\nabla\chi-\epsilon\nabla\chi\otimes\nabla\chi+a\rho^2-\frac{R\rho\theta}{1-b\rho}\big), \\
\displaystyle c_v\big(\rho\theta_t+\rho \mathbf{u}\cdot\nabla\theta\big)+\frac{R\rho\theta}{1-b\rho}\mathrm{div}\mathbf{u}-\kappa(\theta)\Delta\theta=2\nu|\mathbb{D}\mathbf{u}|^2+\lambda(\mathrm{div}\mathbf{u})^2+\mu^2,\\
\displaystyle \rho\chi_t+\rho \mathbf{u}\cdot\nabla\chi=-\mu,\\
\displaystyle \rho\mu=\frac{\rho}{\epsilon}(\chi^3-\chi)-\epsilon \Delta\chi.
\end{array}\right.
\end{equation}
where  $t$ and $ \mathbf{x}\in \mathds{R}^N $ are the time and the spatial variable, $N\geq1$ is  the spatial dimension. $\nabla$, div, and $\Delta$ represent the gradient operator,  the divergence operator and the Laplace operator respectively,  $\cdot_t=\frac{\partial}{\partial t}$. $\rho(\mathbf{x},t)$ and $\theta(\mathbf{x},t)$ are the   density and the  temperature  respectively.  $\mathbf{u}(\mathbf{x},t)$, $\chi(\mathbf{x},t)$ and $\mu(\mathbf{x},t)$ are the velocity, the concentration difference of the immiscible two-phase flow and the chemical potential respectively.   $\nu$, $\lambda$ represent the viscosity coefficients for the immiscible two-phase flow, and $c_v$ is  the specific heat capacity. $\epsilon>0$ is the gradient energy coefficient related to the interfacial width. The stress $p$  is determined through the constitutive \eqref{modified van der Waals equation}, where $v=\frac1\rho$ is the specific volume, and the thermal conductivity coefficient $\kappa$ is of a power-law form with respect to temperature as follows
\begin{equation}\label{kappa}
 \kappa(\theta)=\tilde\kappa\theta^\beta, \qquad\tilde\kappa>0,\qquad\beta\geq0.
\end{equation}
For simplicity, in this paper we assume that $\tilde \kappa=1$.
\begin{remark}
The assumption presented in \eqref{kappa} is theoretically justified by the temperature dependent behavior of thermal conductivity $\kappa$ in compressible immiscible two-phase flows under extreme thermo-fluidic conditions (i.e., high-temperature and high-density regimes). This dependence is rigorously derived from the first-order Chapman-Enskog expansion, which establishes a functional relationship between heat conduction coefficients and temperature. For a comprehensive theoretical foundation, we refer to the seminal work of Chapman-Cowling \cite{CC1994}, where the underlying kinetic theory of gases is systematically developed.
 \end{remark}

\begin{remark}
The so-called immiscible two-phase flow refers to the coexistence of two kinds of immiscible fluids (such as gas-liquid, liquid-liquid) and the existence of a clear phase interface flow, accompanied by the exchange of mass, energy and momentum during the phase transition process. Typical examples include water-steam phase transition, oil-water emulsion flow, etc. The concentration difference  $\chi$ in \eqref{NSFAC3d}  can also be called as the phase field, specifically, taking a micro-volume in the compressible immiscible two-phase flow, assuming that the volume and mass of the micro-volume are $V$ and $M=M_1+M_2$ respectively, where $M_1,M_2$ represent, respectively, the masses of two classes of fluids in the micro-volume of immiscible two-phase flow.   $\chi_i=\frac{M_i}{M}$ $(i=1,2)$ represent the mass concentration of the components in the micro-volume, and the concentration difference of the immiscible two-phase flow is $\chi=\chi_1-\chi_2$. Obviously, physically speaking, $\chi$ should be between $-1$ and $1$, the so-called immiscible two-phase flow diffusion interface, refers to the region  that satisfies $\chi<1$, therefore, we can trace the movement of the diffusion interface through the value of $\chi$.
\end{remark}

It is well known that in recent decades, research on the well-posedness, large-time behavior, and singular limits of the Navier-Stokes equations with the ideal gas state equation has developed rapidly. For instance, the existence of strong solutions under small perturbations was established by Matsumura-Nishida \cite{MN1980}, the existence of renormalized weak solutions was demonstrated by Lions \cite{Lions1998}, and the existence of strong solutions under small perturbations with large oscillations was proven by Huang-Li-Xin \cite{HLX2012}. The global existence of solutions without any restrictions on the smallness of the initial data for the Cauchy problem of one-dimensional compressible non-isentropic Navier-Stokes equations has been successively resolved by Kanel \cite{K1979} and Kazhikhov \cite{K1982}; see Jiang \cite{J1999} and Li-Liang \cite{LL2016} for the large time behavior and further related works. For asymptotic stability results, we refer readers to the studies on rarefaction waves and shock waves by Matsumura-Nishihara \cite{MN1985}-\cite{MN1986}, contact discontinuity waves by Huang-Matsumura-Shi \cite{HMS-2004} and Huang-Li-Matsumura \cite{HLM2010}, interacting shock waves by Huang-Wang-Wang-Yang \cite{HWWY2015}, and the vanishing viscosity limit for interacting shock waves by Shi-Yong-Zhang \cite{SYZ2016}, among others.
  
The research on \eqref{NSFAC3d} with the ideal gas equation of state, primarily focuses on the existence theory of solutions. For instance, Feireisl-Petzeltová-Rocca-Schimperna \cite{FPRS2010} proved the global existence of renormalized finite-energy weak solutions under the condition that the adiabatic exponent of pressure satisfies $\gamma > 6$. This result was subsequently extended to $\gamma > 2$ by Chen-Wen-Zhu \cite{cwz2019}. Additionally, Kotschote \cite{Kotschote2012} and Chen-He-Huang-Shi \cite{CHHS2022-1} established the local strong solution existence theory. In one-dimensional cases, Ding-Li-Luo \cite{DLL2013} investigated the global existence of strong solutions for the Cauchy problem, while Ding-Li-Tang \cite{dlt2019} addressed the free boundary problem. The initial density containing vacuum was studied by Chen-Guo \cite{CG2017} and Li-Yan-Ding-Chen \cite{lydc2023}, and the dependence of heat conduction coefficients on temperature was explored by Chen-He-Huang-Shi \cite{CHHS2021, CHHS2022}. Furthermore, Yan-Ding-Li \cite{ydl2022} examined the case where viscosity depends on the phase field. For the three-dimensional Cauchy problem with small initial perturbations, Chen-Hong-Shi \cite{CHS-2021, CHS-2023}, Chen-Li-Tang \cite{clt2021-non-nsac}, and Zhao \cite{zhao2019} analyzed the global existence and uniqueness of strong solutions. Lastly, the large-time behavior of rarefaction waves, contact waves, or stationary solutions under specific restrictions on the initial phase field was considered by Yin-Zhu \cite{yz2019}, Luo-Yin-Zhu \cite{lyz2018, lyz2020}, and Chen-Peng-Shi \cite{CPS-2024}, among others.

However, for fluids undergoing phase transitions, the pressure $p$ is typically required to be non-monotonic with respect to  $v$, referred to as the non-ideal fluid equation of state. For the isentropic compressible Navier-Stokes system with van der Waals equation of state, Hsieh-Wang \cite{HW-1997} conducted a numerical analysis using a pseudo-spectral method with artificial viscosity. Their results demonstrated that the non-monotonicity of pressure induces phase transitions. He-Liu-Shi \cite{HLS} replicated Hsieh-Wang's work by employing a second-order TVD Runge-Kutta splitting scheme combined with the Jin-Xin relaxation scheme, which was later extended to the compressible Navier-Stokes/Cahn-Hilliard system by He-Shi \cite{HS-2020}. For dynamic elastic bar theories where the stress-deformation relation is non-monotone, such as the pressure $p(v)=v^{3}-v$, Mei-Liu-Wong \cite{MLW-2007-1}-\cite{MLW-2007-2} established the existence of global strong solutions for the one-dimensional Navier-Stokes system with additional artificial viscosity. Hoff-Khodia \cite{HK-1993} investigated the dynamic stability of certain steady-state weak solutions in one-dimensional whole space under small initial perturbations.

It is well established that phase transitions give rise to novel interfacial structures, yet the underlying mechanisms governing their formation and the nature of the resultant interfaces, particularly whether they retain diffusion interface characteristics remain poorly understood. For another mathematical model of non-homogeneous two-phase flow, the Navier-Stokes/Cahn-Hilliard system, when the state equation is that of an ideal fluid, He-Shi \cite{HS-2021} coupled a fully discrete local discontinuous Galerkin  finite element method with a scalar auxiliary variable  approach to numerically resolve the compressible Navier-Stokes/Allen-Cahn system under van der Waals equation of state. Their computational framework demonstrated that initial thermodynamic conditions critically dictate the long-term evolution of diffusive interfaces in such systems, offering valuable insights into phase behavior under compressibility constraints. Notably, investigations into compressible immiscible two-phase flows governed by non-ideal fluid equations of state remain sparse.

 The motivation of  this paper is to understand the phase transition in immiscible two-phase flow and its effect on the diffusion interface.
Further, by introducing the Lagrange coordinate system (without causing confusion, we still write the new coordinate system as $x,t$, and the viscosity coefficient is assumed to be $1$),   the Cauchy problem for system \eqref{NSFAC3d} in one dimension can be rewritten as
\begin{equation}\label{NSFAC-Lagrange}
\left\{\begin{array}{llll}
\displaystyle v_t-u_x=0,\\
\displaystyle u_t+\big(\frac{R\theta}{v-b}-\frac a{v^2}\big)_x=\big(\frac{u_{x}}{v}\big)_x-\frac\epsilon2\big(\frac{\chi_x^2}{v^2}\big)_x, \\
\displaystyle \theta_t+\frac{R\theta}{v-b}u_x-\big(\frac{\kappa(\theta)\theta_x}{v}\big)_x=\frac{u_x^2}{v}+v\mu^2,\\
\displaystyle \chi_t=-v\mu,\\
\displaystyle \mu=\frac{1}{\epsilon}(\chi^3-\chi)-\epsilon\big(\frac{\chi_x}{v}\big)_x,
\end{array}\right.
\end{equation}
with initial condition
\begin{equation}\label{initial condition}
 (v,u,\theta,\chi)(x,0)=(v_0,u_0,\theta_0,\chi_0)(x)\xrightarrow{x\rightarrow\pm\infty}(\bar{v},0,\bar{\theta},\pm1).
\end{equation}
Here $\bar v$ and $\bar \theta$ are two positive constants satisfying one of the following constraints
\begin{equation}\label{bar v, bar theta-1}
  \frac{8a}{27Rb}>\bar\theta>0,\qquad  b+h<\bar v<v_*(\bar\theta), \ \mathrm{or} \ \bar v>v^*(\bar\theta),
\end{equation}
or
\begin{equation}\label{bar v, bar theta-2}
  \bar\theta\geq\frac{8a}{27Rb},\qquad \bar v>b+h,
\end{equation}
where  $v_*(\bar\theta)$ and $v^*(\bar\theta)$ in \eqref{bar v, bar theta-1} are given in Lemma \ref{properities of p}.

 \begin{remark}When the temperature is lower than the critical temperature (i.e., $\theta<\frac{8a}{27Rb}$, see Figure 1.(a)), the van der Waals equation of state \eqref{modified van der Waals equation} is non-monotonic with respect to volume $v$. At this time, if the constraint on $\bar v$ in \eqref{bar v, bar theta-1} does not hold, the positivity of \eqref{Phi} cannot be guaranteed, and thus the global existence of the solution cannot be ensured. When the temperature surpasses the critical temperature (i.e., $\theta\geq\frac{8a}{27Rb}$, see Figure 1.(a)),  the van der Waals equation of state exhibits a distinct monotonic relationship with respect to the volume $v$. In this supercritical state, no matter which value $\bar v>b+h$ is chosen, the \eqref{Phi} remains positive. This classification results in three regions for initial states regions $D_{unstable},D_{metastable},D_{stable}$  corresponding to unstable, metastable, and stable regions(see depicted in Fig. 1(b)).
\begin{equation}\label{ABC}
\left.\begin{array}{llll}
\displaystyle D_{unstable}=\big\{v\big|v_\alpha(\bar\theta)<v<v_\beta(\bar\theta)\big\} & \mathrm{unstable \ region},\\
\displaystyle D_{metastable}=\big\{v\big|v_*(\bar\theta)<v<v_\alpha(\bar\theta),v_\beta(\bar\theta)<v<v^*(\bar\theta)\big\},& \mathrm{metastable \ region},\\
\displaystyle  D_{stable}=\big\{v\big|b+h<v<v_*(\bar\theta)<v_\alpha(\bar\theta),v>v^*(\bar\theta)\big\},& \mathrm{stable\ region}.
\end{array}\right.
\end{equation}
In summary, the classification \eqref{ABC} indicates that the initial $v_0$ at infinity must lie within the stable region $D_{stable}$, which a fundamental requirement in our framework. Moreover,
existing computational work (see Hsieh-Wang \cite{HW-1997}, He-Shi \cite{HS-2021}) indicates that the density of the initial state determines the long-term behavior of the system also. A more detailed discussion on this classification will be provided in our subsequent work.
 \end{remark}

\

\noindent\textbf{\normalsize Notations.} Throughout our paper, we denote  $L^\infty(\mathbb{R})$ is the space of all bounded measurable functions on $\mathbb{R}$, with the norm $\|f\|_{L^\infty}$ = esssup$_{x\in\mathbb{R}}|f|$,  $L^q(\mathbb{R})$ ($q\geq1$, and $L^q$ without any ambiguity)  as the space of $q$-times integrable real valued function defined on $\mathbb{R}$, with the norm $\|f\|_{L^q}=(\int_{\mathbb{R}}|f|^qdx)^{\frac{1}{q}}$, For the special case that $q=2$, we simply denote $\|f\|_{L^2}$ by $\|f\|$. We further denote $H^l(\mathbb{R})$ ($l>0, H^l$ without any ambiguity)  the Sobolev space of $L^2$-functions $f$ on $\mathbb{R}$ whose derivatives $\partial^j_x f,  j=1,\cdots,l$ are also square integrable functions, with the norm $ \|f\|_{l}=(\sum_{j=0}^l\|{\partial^j_x f}\|^2)^{\frac{1}{2}}$. In addition, we denote
by $L^\infty([0, T]; H^k(\mathbb{R}))$ (resp. $L^2(0, T; H^k(\mathbb{R})))$ the space of bounded measurable (resp. square integrable) functions on $[0,T]$ with values taken in a Banach space $H^k(\mathbb{R})$.

\ 

In this paper, we adopt the following fundamental assumptions concerning the initial conditions of the system \eqref{NSFAC-Lagrange} under investigation
\begin{equation}\label{condition 1}v_0-\bar v\in L^1(\mathbb{R}),\quad
(v_0-\bar v,u_0,\chi_{0x})\in H^1(\mathbb{R}),\quad \theta_0-\bar \theta\in H^1(\mathbb{R}),\quad \chi^2_0-1\in L^2(\mathbb{R}),
\end{equation}
and
\begin{equation}\label{condition 2}
 \inf_{x\in\mathbb{R}}v_0(x)>b+h,\ \ \inf_{x\in\mathbb{R}}\theta_0(x)>0,\ \ \chi_0(x)\in[-1,1],\quad \forall x\in\mathbb{R}.
\end{equation}

 \

Now, we can present the  following result on the existence and uniqueness of the strong solution for the Navier-Stokes/Allen-Cahn system \eqref{NSFAC-Lagrange}-\eqref{initial condition} with van der Waals equation of state \eqref{modified van der Waals equation} and  degenerate thermal conductivity $\kappa(\theta)=\theta^\beta, (\beta>0)$ .
\begin{theorem}\label{thm-global} Under the assumptions of \eqref{modified van der Waals equation}--\eqref{Hypothesis of p}, \eqref{bar v, bar theta-1}--\eqref{bar v, bar theta-2}, and \eqref{condition 1}--\eqref{condition 2}, the thermal conductivity coefficient $\kappa(\theta)=\theta^\beta,(\beta>0)$,
then the Cauchy problem of the compressible non-isentropic Navier-Stokes/Allen-Cahn system \eqref{NSFAC-Lagrange}-\eqref{initial condition} has a unique global strong solution $(v,u,\theta,\chi)$ such that,  for any $T>0$, there exists a constant $C>0$ such that
\begin{equation}
\left.\begin{array}{llll}
 \displaystyle v-\bar v\in L^1(\mathbb{R}),\ \big(v-\bar v,u,\theta-\bar\theta,\chi_x\big)\in L^\infty(0,T;H^1(\mathbb{R})),\  \chi^2-1\in L^\infty(0,T;L^2(\mathbb{R})), \\
 \displaystyle  v_x\in  L^2\big(0,T;L^2(\mathbb{R})\big), \quad (u_x,\theta_x,\chi_x)\in  L^2\big(0,T;H^1(\mathbb{R})\big),\quad \big(\frac {\chi_{x}}{v}\big)_{xx}\in  L^2\big(0,T;L^2(\mathbb{R})\big),  \\
 \displaystyle (v_t,\chi_t)\in L^\infty(0,T;L^2(\mathbb{R}))\cap L^2(0,T;H^1(\mathbb{R})), \quad\big(u_t,\theta_t,\chi_{xt}\big)\in L^2(\mathbb{R}\times(0,T)),
\end{array}\right.
\end{equation}
and
\begin{equation}\label{upper and lower bound}
b+h< v(x,t)\leq C,\ \ C^{-1}\leq \theta(x,t)\leq C,\ \ \chi(x,t)\in[-1,1],\quad (x,t)\in\mathbb{R}\times[0,T],
\end{equation} where the positive constant $C$ only depends on $T$, the initial data $v_0,u_0,\theta_0,\chi_0$ and $\epsilon,h$.
\end{theorem}

\

The compressible Navier-Stokes equations emerge as a fundamental special case of system \eqref{NSFAC-Lagrange} when modeling single-fluid flows, corresponding to the parameter choice $\chi\equiv1$  or $\chi\equiv-1$. Under the Lagrange coordinate system, the one dimensional Cauchy problem of compressible Navier-Stokes system governed by the van der Waals equation of state \eqref{modified van der Waals equation} and degenerate thermal conductivity $\kappa(\theta)=\theta^\beta, (\beta\geq0)$ is as following
\begin{equation}\label{NSF-Lagrange}
\left\{\begin{array}{llll}
\displaystyle v_t-u_x=0,\\
\displaystyle u_t+\big(\frac{R\theta}{v-b}-\frac a{v^2}\big)_x=\big(\frac{u_{x}}{v}\big)_x, \\
\displaystyle \theta_t+\frac{R\theta}{v-b}u_x-\big(\frac{\kappa(\theta)\theta_x}{v}\big)_x=\frac{u_x^2}{v},
\end{array}\right.
\end{equation}
with initial condition
\begin{equation}\label{initial condition of NSF}
 (v,u,\theta)(x,0)=(v_0,u_0,\theta_0)(x)\xrightarrow{x\rightarrow\pm\infty}(\bar{v},0,\bar{\theta}).
\end{equation}
By adjusting the proof framework in Theorem \ref{thm-global}, we can obtain the following existence and uniqueness result for the strong solution of the Navier-Stokes system \eqref{NSF-Lagrange}-\eqref{initial condition of NSF} governed by the van der Waals equation of state \eqref{modified van der Waals equation} and degenerate thermal conductivity $\kappa(\theta)=\theta^\beta, (\beta\geq0)$.
\begin{theorem}\label{NS} Under the assumptions of \eqref{modified van der Waals equation}--\eqref{Hypothesis of p}, \eqref{bar v, bar theta-1}-\eqref{bar v, bar theta-2}, the thermal conductivity coefficient  $\kappa(\theta)=\theta^\beta,(\beta\geq0)$, and
\begin{equation}\label{condition 1'}
\left.\begin{array}{llll}
\displaystyle v_0-\bar v\in L^1(\mathbb{R}),\
v_0-\bar v,u_0\in H^1(\mathbb{R}),\quad \theta_0-\bar \theta\in H^1(\mathbb{R}),\\
\displaystyle \inf_{x\in\mathbb{R}}v(x)>b+h,\ \ \inf_{x\in\mathbb{R}}\theta_0(x)>0,
\end{array}\right.
\end{equation}
then the Cauchy problem of the compressible non-isentropic Navier-Stokes system \eqref{NSF-Lagrange}-\eqref{initial condition of NSF} has a unique global strong solution $(v,u,\theta)$ such that,  for any $T>0$, there exists a constant $C>0$ such that
\begin{equation}
\left.\begin{array}{llll}
 \displaystyle v-\bar v\in L^{\infty}(0,T;L^1(\mathbb{R})),\ \ \big(v-\bar v,u,\theta-\bar\theta\big)\in L^\infty(0,T;H^1(\mathbb{R})), \\
 \displaystyle  v_x\in  L^2\big(0,T;L^2(\mathbb{R})\big), \quad (u_x,\theta_x)\in  L^2\big(0,T;H^1(\mathbb{R})\big),  \\
 \displaystyle v_t\in L^\infty(0,T;L^2(\mathbb{R}))\cap L^2(0,T;H^1(\mathbb{R})), \quad\big(u_t,\theta_t\big)\in L^2(\mathbb{R}\times(0,T)).
\end{array}\right.
\end{equation}
Moreover
\begin{equation}\label{upper and lower bound}
b+h\leq v(x,t)\leq C,\ \ 0<C^{-1}\leq \theta(x,t)\leq C,\quad (x,t)\in\mathbb{R}\times[0,T],
\end{equation} where the positive constant $C$ only depends on $T$ , the initial data $v_0,u_0,\theta_0$ and $h$.
\end{theorem}

\begin{remark} Theorem \ref{thm-global} provides  the global existence and uniqueness of strong solutions to the one-dimensional Cauchy problem of Nqvier-Stokes/Allen-Cahn system \eqref{NSFAC-Lagrange}-\eqref{initial condition} for non-vacuum initial data without requiring smallness conditions on the initial conditions. Specifically, the following key results hold:
\begin{enumerate}
  \item[(i)] Despite the non-monotonic pressure inducing significant density variations, triggering a phase transition in the immiscible two-phase flow, the density and temperature remain strictly bounded away from zero and infinity for all finite times. Crucially, this guarantees the absence of vacuum formation.

 \item[(ii)] Owing to the stabilizing effects of physical viscosity, both the density and phase field retain their continuity even in proximity to the phase transition region. This regularity persists despite the sharp gradients induced by the phase transition.

 \item[(iii)] A direct corollary is that any phase transition within the compressible immiscible two-phase system necessitates the emergence of a new diffusion interface at the transition locus. This interface arises as a dynamical consequence of the interplay between pressure non-monotonicity and viscous dissipation.

\item[(iv)] Within any finite time interval, the phase field of the compressible immiscible two-phase flow remains continuous, even when the diffusion interfaces interact. Furthermore, density, temperature, phase field and velocity of the immiscible two-phase system  all maintain continuity throughout this process. 
\end{enumerate}
 \end{remark}

\begin{remark} To rigorously control the density bounds, particularly the oscillatory behavior of the density gradient in phase transition regions, prior studies such as Hsieh-Wang \cite{HW-1997} introduced an artificial viscous term involving the second derivative of density. This term was essential for stabilizing numerical solutions of the compressible Navier-Stokes system. A similar approach was later adopted by Mei-Liu-Wong \cite{MLW-2007-1,MLW-2007-2}, where the inclusion of such artificial dissipation became unavoidable in their analysis.
However, this artificial viscosity lacks physical justification, as it introduces non-physical damping effects that may distort the true dynamics of phase transitions. In contrast, our work presents a novel analytical framework that circumvents this limitation. Through a combination of refined energy estimates and carefully designed functional spaces, we achieve robust control over density oscillations without relying on any unphysical regularization terms. This advancement not only preserves the mathematical consistency of the system but also ensures that the solutions remain faithful to the underlying physics.
\end{remark}

\begin{remark} The condition that $\beta>0$ in Theorem \ref{thm-global}  is necessary in this paper.
 The reason is that ,
 the interface free energy $\frac{\epsilon}{2}\big(\frac{\chi_{x}^2}{v^2}\big)_{x}$ that appears in the momentum equation \eqref{NSFAC-Lagrange}$_2$ requires a higher-order estimate of the phase field $\chi$, but the phase field $\chi$ is also coupled with $v$ in \eqref{NSFAC-Lagrange}$_5$, it is necessary to obtain higher-order a priori estimates of $v$ simultaneously. However, the mass conservation equation \eqref{NSFAC-Lagrange}$_1$ for $v$ is a hyperbolic equation, and the regularity of the solution $v$ itself is weak. Therefore, these difficulties lead to the situation that applying the techniques presented in this paper under the conditions of Theorem \ref{thm-global} cannot yield an upper bound for $v$.
   However, for single-phase compressible fluids without a diffusion interface, since there is no influence of the phase field $\chi$,  mathematically, it is manifested that the Naver-Stokes/Allen-Cahn system \eqref{NSFAC-Lagrange} degenerates into the Navier-Stokes system \eqref{NSF-Lagrange}, the aforementioned problem no longer exists. This is the fundamental reason why Theorem \ref{NS} can be derived by adjusting the proof framework.
\end{remark}

\begin{remark}
Theorem \ref{NS} establishes a rigorous mathematical framework to validate the well-posedness and physical consistency of the global strong solutions derived by Hsieh-Wang \cite{HW-1997}, He-Shi \cite{HS-2021}. However, the long-time behavior of these solutions remains beyond the scope of this paper, as it involves the Riemann problem for mixed-type systems of conservation laws, a topic of significant complexity and theoretical challenge. A comprehensive analysis of this issue will be addressed in future work.
\end{remark}

At the end of this section, we will briefly  emphasize the key steps in the proof process, and outline the main steps of the proof presented in this paper.

The key to the proofs of Theorem 1.1 and Theorem 1.2 lay in providing an upper bound for the volume $v$ and the upper and lower bounds for the temperature $\theta$.
For this purpose, we successfully observed and constructed two energy functionals $\Phi(v)$ and $\Psi(\theta)$ (see \eqref{Phi}-\eqref{Psi}), from which the standard energetic estimate \eqref{basic energy inequality} was derived. For the immiscible two-phase flow discussed in Theorem 1.1, we only consider the case where the exponent $\beta>0$ of the thermal conductivity coefficient $\kappa=\theta^\beta$. At this point, the boundedness of $\|\frac{\theta^\beta\theta_x^2}{v\theta^2}\|_{L^1(0,T;L^1(\mathbb{R}))}$   can be obtained from the  standard energetic estimate \eqref{basic energy inequality}, therefore, the truncation technique introduced by Jiang \cite{J1999} and Li-Liang \cite{LL2016} (see\eqref{bar v theta} and \eqref{the truncation function method}) can thus be adopted to obtain the lower bound of temperature $\theta$, and in particular, two key energy inequalities \eqref{integral sup theta} and \eqref{theta0} were obtained as a result. Furthermore, based on the expression of $v$ \eqref{expression of v} that we derived with the inspiration from  Kazhikhov \cite{K1982},  the  upper bound of $v$ can be obtained. 
For the case where $\beta= 0$, due to the reasons stated in Remark 1.6, at present, we can only derive the well-posedness analysis for the compressible non-isentropic Navier-Stokes equations that describe single-phase flow. Even so,
 since there is no longer a boundedness property for $\|\frac{\theta^\beta\theta_x^2}{v\theta^2}\|_{L^1(0,T;L^1(\mathbb{R}))}(\beta>0)$ for this case, this led to the failure of \eqref{eta} and \eqref{The integral of theta}. Therefore, inspired by the arguments from  L$\mathrm{\ddot{u}}$-Shi-Xiong \cite{lsx2021}, we have developed a new complex but effective energy estimation procedure to obtain the upper bound for $v$.
Moreover,
To obtain the higher-order estimates of $v,u, \chi,\theta$, taking into account the characteristics of van der Waals fluids, we employ the pressure decomposition technique (see \eqref{Another form of the momentum equation}), combined with energy estimation, the boundedness of $\|v_x\|_{L^2(0,T;L^2(\mathbb{R}))}$ is achieved. Furthermore, since the interface free energy $\frac{\epsilon}{2}\big(\frac{\chi_{x}^2}{v^2}\big)_{x}$ appears in the momentum equation \eqref{NSFAC-Lagrange}$_2$, just the fact that only the gradient estimation of the phase field $\chi$  is not sufficient to guarantee the attainment of the upper bound of temperature $\theta$. Therefore, we  transform \eqref{NSFAC-Lagrange}$_5$ into $\frac{\chi_t}{v}-\epsilon\big(\frac{\chi_x}{v}\big)_{x}=-\frac{1}{\epsilon}(\chi^3-\chi)$,  and through the elementary but technical $L^2$ energy method,  $\big(\frac{\chi_{x}}{v}\big)_{xx}\in L^2(0,T;L^2(\mathbb{R}))$ was obtained. Moreover, by continuing to apply the truncation technique and combining with energy estimation, the boundedness of $\|u_x\|_{L^2(0,T;L^2(\mathbb{R}))}$ was derived, and finally, multiplying \eqref{NSFAC-Lagrange}$_3$ by $\theta^\beta\theta_t$,  one can derive the upper bound of $\theta$.

\
   
The remainder of this paper is organized as follows. In Section \ref{sec2}, the local existence and uniqueness of solutions for  Navier-Stokes/Allen-Cahn system \eqref{NSFAC-Lagrange}-\eqref{initial condition} is established via the  linearization method and iteration technique. In Section \ref{sec3}, we consider the case where the thermal conductivity coefficient $\kappa(\theta)=\theta^\beta, (\beta>0)$. Building upon the local solution established in Section \ref{sec2}, we derive a priori estimates that only depend on $T$, the initial data $v_0,u_0,\theta_0,\chi_0$ and $\epsilon,h$. These estimates ultimately enable us to obtain the existence and uniqueness of global solutions to the Navier-Stokes/Allen-Cahn system \eqref{NSFAC-Lagrange}-\eqref{initial condition}. In Section \ref{sec4}, when the thermal conductivity $\kappa$ is a positive constant,  the analytical framework in  Section \ref{sec3} no longer applies to obtaining the upper bound for $v$. We therefore develop an alternative energy estimation method (more intricate but effective) to overcome this difficulty, and thus,  the existence and uniqueness of global-in-time solution is obtained for  compressible Navier-Stokes system  \eqref{NSF-Lagrange}-\eqref{initial condition of NSF}.

\section{Existence of Local Strong Solutions}\label{sec2}
\ \ \ \ In this section, we prove the existence and uniqueness of local strong solutions to the Cauchy problem for the system \eqref{NSFAC-Lagrange}-\eqref{initial condition} via an iterative method. Notably, the Navier-Stokes system \eqref{NSF-Lagrange}-\eqref{initial condition of NSF} is a special case of the Navier-Stokes/Allen-Cahn system \eqref{NSFAC-Lagrange}-\eqref{initial condition}; thus, the local well-posedness result directly applies to the former. For brevity, we omit this straightforward extension in later analysis. Prior to the detailed proof, we define the following solution space
 $X_{m_1,m_2,M}(0,T)$  as follows
\begin{eqnarray}\label{Solution space}
 &X_{m_1,m_2,M}(0,T)=\displaystyle\Big\{(v, u,\theta, \chi)\:\Big| (v-\bar v,u,\theta-\bar\theta,\chi_x) \in C\big(0,T;H^1(\mathbb{R})\big),\notag\\
 &\qquad\displaystyle\quad\qquad\qquad \chi^2-1\in C\big(0,T;L^{2}(\mathbb{R})\big),\ \inf_{\mathbb{R}\times[0,T]}\theta\geq m_1>0,\ \inf_{\mathbb{R}\times[0,T]}v\geq m_2>b+h>0,\notag\\
&\qquad\qquad\qquad\displaystyle\sup_{ t \in [0,T]}\big(\|(v-\bar v,u,\theta-\bar\theta,\chi_x)\|_1+\|\chi^2-1\|\big)\leq M,\\
&\qquad\qquad\qquad\displaystyle\qquad\qquad\quad \big(v-\bar v,u,\theta - \bar\theta\big)_x\in L^2\big(0,T;H^1(\mathbb{R})\big),\chi_x\in L^2\big(0,T;H^2(\mathbb{R})\big)\Big\}.\notag
\end{eqnarray}
for any given  $M>0,m_1>0,m_2>b+h$.

Now, through the linearization method and iteration technique, we will establish the existence and  uniqueness of local strong solutions for the system \eqref{NSFAC-Lagrange}-\eqref{initial condition} as following.

\begin{proposition}  \label{local existence}
	\textbf{(Local Existence).}  Under the assumptions of \eqref{modified van der Waals equation}--\eqref{Hypothesis of p}, \eqref{bar v, bar theta-1}-\eqref{bar v, bar theta-2}, $\kappa=\theta^\beta,\beta\geq0$, and
 \begin{equation}\label{condition-11 for local solution}
  \| (v_0-\bar v,u_{0},\theta_{0}-\bar\theta,\chi_{0x})\|_1+\|\chi_0^2-1\|\leq M,
 \end{equation}
 and
  \begin{equation}\label{condition-2 for local solution}
   \inf\limits_{x\in \mathbb{R} }\theta_0(x)\geq m_1>0,\ \inf\limits_{x\in \mathbb{R} }v_0(x)\geq m_2>b+h,
  \end{equation}
then there exists a small time  $T^*$ depending only on $m_1,m_2,M$, and $v_0,u_0,\theta_0,\chi_0$ such that, the Cauchy problem \eqref{NSFAC-Lagrange}-\eqref{initial condition}, admits a unique solution $(v,u, \theta, \chi)\in X_{\frac{m_1}{2},\frac{m_2+b+h}{2},2M}\big(0,T^*\big)$.
\end{proposition}

\begin{proof}  
From the system \eqref{NSFAC-Lagrange}$_1$, one has
\begin{equation}\label{Approximate equation 1}
  v(x,t)=v_0(x)+\int_0^tu_x(x,\tau)d\tau.
\end{equation}
Rewriting  \eqref{NSFAC-Lagrange}$_2$, \eqref{NSFAC-Lagrange}$_3$ and \eqref{NSFAC-Lagrange}$_{4,5}$  as the following   equations of $u$, $\theta$ and $\chi$ respectively
\begin{equation}\label{Approximate equation 2}
\left\{\begin{array}{llll}
\displaystyle  u_t-\big(\frac{u_x}{v}\big)_x=g_1(v,v_x,\theta,\theta_x,\chi_x,\chi_{xx}),    \\
\displaystyle u(x,0)=u_0(x),
\end{array}\right.
\end{equation}
\begin{equation}\label{Approximate equation 3}
\left\{\begin{array}{llll}
\displaystyle \theta_t-\big(\frac{\theta^\beta\theta_x}{v}\big)_x=g_2(v,\theta,v_x,u_x,\chi,\chi_x,\chi_{xx}),\\
\displaystyle \theta(x,0)=\theta_0(x),
\end{array}\right.
\end{equation}
and
\begin{equation}\label{Approximate equation 4}
\left\{\begin{array}{llll}
\displaystyle \chi_t-\epsilon\chi_{xx}=g_3(v,v_x,\chi,\chi_x),\\
\displaystyle\chi(x,0)=\chi_0(x),
\end{array}\right.
\end{equation}
where $g_1$, $g_2$, $g_3$ are 
\begin{equation}\label{g1}
  g_1\xlongequal{\mathrm{def}} -\frac{\epsilon}{2}\big(\frac{\chi_x^2}{v^2}\big)_x-\big(-\frac{a}{v^2}+\frac{R\theta}{v-b}\big)_x,
\end{equation}
and
\begin{equation}
\label{g2 and g3}
g_2\xlongequal{\mathrm{def}} v\big[\frac1\epsilon(\chi^3-\chi)-\epsilon(\frac{\chi_x}{v})_x\big]^2+\frac{u_x^2}{v}-\frac{R\theta}{v-b}u_x,\qquad
 g_3\xlongequal{\mathrm{def}} -\frac{v}{\epsilon}(\chi^3-\chi)-\frac{\epsilon\chi_x v_x}{v}.
\end{equation}

In order to  construct a Cauchy sequence for the approximate solutions of \eqref{Approximate equation 1}-\eqref{Approximate equation 4},  we plan to use the iterative method, for which we construct approximate sequence $\{(v_{0k},u_{0k},\theta_{0k},\chi_{0k})\}_{k=1}^{\infty}$ for initial data  $(v_0-\bar v,u_0,\theta_0-\bar\theta)\in H^1(\mathbb{R})$, $\chi_{0x}\in H^1(\mathbb{R})$ and $\chi_{0}^2-1\in L^2(\mathbb{R})$, such that 
\begin{equation}\label{approximate sequences of initial data}
  (v_{0k}-\bar v,u_{0k},\theta_{0k}-\bar\theta)\in H^4(\mathbb{R}),\ (\chi_{0k})_x\in H^3(\mathbb{R}),\ (\chi_{0k})^2-1\in L^2(\mathbb{R}),
\end{equation}
satisfying 
\begin{equation}\label{initial approximate sequence}\left.\begin{array}{llll}
\displaystyle \|(v_{0k}-\bar v,u_{0k},\theta_{0k}-\bar\theta)\|_2+\|\chi_{0k}^2-1\|+\| (\chi_{0k})_x\|_1\leq \frac{3}{2}M,\\
\displaystyle \inf_{x\in \mathbb{R} }\theta_{0k}\geq\frac{2}{3}m_1,\quad\inf_{x\in \mathbb{R} }v_{0k}\geq\frac{m_2}{3}+\frac{2(b+h)}{3},\ \ \forall k=1,2,\cdots,\end{array}\right.
\end{equation}
and 
\begin{equation}\label{convergence}\left.\begin{array}{llll}
 (v_{0k},u_{0k},\theta_{0k})\xrightarrow{k\rightarrow\infty}(v_0,u_{0},\theta_{0}),\ \mathrm{strongly\ in}\ H^1(\mathbb{R}),\\
  (\chi_{0k})_x\xrightarrow{k\rightarrow\infty}(\chi_{0})_x,\ \mathrm{strongly\ in}\ H^1(\mathbb{R}),\ \chi_{0k}^2-1\xrightarrow{k\rightarrow\infty}\chi_{0}^2-1,\ \mathrm{strongly\ in}\ L^2(\mathbb{R}).
\end{array}\right.\end{equation}
For fixed $k=1,2,3,\cdots$, letting the corresponding sequence of approximate solutions 
\begin{equation}\label{sequence of approximate solutions}\left.\begin{array}{llll}
\displaystyle\big\{(v_k^{(n)},u_k^{(n)},\theta_k^{(n)},\chi_k^{(n)})(x,t)\big\}_{n=1}^{+\infty},\qquad\big(v_k^{(0)},u_k^{(0)},\theta_k^{(0)},\chi_k^{(0)}\big)=(v_{0k},u_{0k},\theta_{0k},\chi_{0k}),
\end{array}\right.\end{equation}
where $u_k^{(n)}$, $\theta_k^{(n)}$ and $\chi_k^{(n)}$ are given by the linearized iterative scheme for \eqref{Approximate equation 1}-\eqref{Approximate equation 4} as following
\begin{equation}\label{iterative format 2}
\left\{\begin{array}{llll}
\displaystyle  u_{kt}^{(n)}-\Big(\frac{u^{(n)}_{kx}}{v_k^{(n-1)}}\Big)_x=g_1(v_k^{(n-1)},v_{kx}^{(n-1)},\theta_k^{(n-1)},\theta_{kx}^{(n-1)},\chi_{kx}^{(n-1)},\chi_{kxx}^{(n-1)}),    \\
\displaystyle u_k^{(n)}(x,0)=u_{0k}(x),
\end{array}\right.
\end{equation}

\begin{equation}\label{iterative format 3}
\left\{\begin{array}{llll}
\displaystyle \theta^{(n)}_{kt}-\Big(\frac{(\theta^{(n-1)})^\beta\theta^{(n)}_{kx}}{v^{(n-1)}_k}\Big)_x =g_2(v_k^{(n-1)},\theta_k^{(n-1)},v^{(n-1)}_{kx},u_{kx}^{(n-1)},\chi_{k}^{(n-1)},\chi_{kx}^{(n-1)},\chi_{kxx}^{(n-1)}),\\
\displaystyle \theta_k^{(n)}(x,0)=\theta_{0k}(x),
\end{array}\right.
\end{equation}

\begin{equation}\label{iterative format 4}
\left\{\begin{array}{llll}
\displaystyle \chi^{(n)}_{kt}-\chi_{kxx}^{(n)}=g_3(v_k^{(n-1)},v_{kx}^{(n-1)},\chi_{k}^{(n-1)},\chi_{kx}^{(n-1)}),\\
\displaystyle\chi^{(n)}_k(x,0)=\chi_{0k}(x),
\end{array}\right.
\end{equation}
and $ v^{(n)}_{k}$ is given as following
\begin{equation}\label{iterative format 1}
   v^{(n)}_{k}(x,t)=v_{0k}(x)+\int_0^tu^{(n)}_{kx}(x,\tau)d\tau.
\end{equation}
Using the theory of second order linear parabolic equation,  if
    \begin{equation}\label{g1g2g3}
    g_i^{(n-1)}\in C([0,T];H^2(\mathbb{R})), \ \ \ i=1,2,3,
  \end{equation}
and
\begin{equation}\label{First iteration}
 (v_{0k}-\bar v,u_{0k},\theta_{0k}-\bar\theta)\in H^4(\mathbb{R}),\  (\chi_{0k})_x\in H^3(\mathbb{R}),\  (\chi_{0k})^2-1\in L^2(\mathbb{R}),
\end{equation}
 then, combining with the maximum principle for parabolic equation (see Lemma 2.1 in \cite{p2005}), the equations \eqref{iterative format 2}-\eqref{iterative format 4} admit a unique local solution $(u^{(n)}_k,\theta^{(n)}_k,\chi_k^{(n)})$, such that
\begin{eqnarray}\label{local solution}
&& (u^{(n)}_k,\theta^{(n)}_k-\bar\theta)\in C([0,T];H^4(\mathbb{R}))\cap C^1([0,T];H^2(\mathbb{R}))\cap L^2([0,T];H^5(\mathbb{R})),\notag\\
&& \chi_{kx}^{(n)}\in C([0,T]; H^3(\mathbb{R}))\cap L^2([0,T];H^4(\mathbb{R})),\ \ \chi_{kt}^{(n)}\in C^0([0,T];H^2(\mathbb{R})),\notag\\
&&  (\chi_k^{(n)})^2-1\in C([0,T]; H^2(\mathbb{R})). \notag
\end{eqnarray}
Further, if $(v_k^{(n-1)},u_k^{(n-1)},\theta_k^{(n-1)},\chi_k^{(n-1)})\in X_{\frac{1}{2}m_1,\frac{m_2+b+h}{2},2M}(0,t_0)$, from the elementary energy estimate, choosing $t_0=t_0(m,M)$  suitably small,  one has
\begin{equation}\label{local estimate-1}
  \|u_k^{(n)}(t)\|_2^2\leq \big((\frac{3}{2}M)^2+C(m,M)t_0\big)e^{C(m,M)t_0}\leq(2M)^2,
\end{equation}
and 
\begin{equation}\label{local estimate-2}
  \int_0^{t_0}\|u_{kx}^{(n)}(\tau)\|_2^2d\tau\leq C(m,M)(2M)^2,
\end{equation}
similarly, one gets
\begin{equation}\label{local estimate-3}
  \big\|\theta_k^{(n)}-\bar\theta\big\|_2^2\leq (2M)^2,\ \ \big\|(\chi_{k}^{(n)})^2-1\big\|^2\leq (2M)^2, \ \  \big\|\chi_{kx}^{(n)}\big\|_2^2\leq (2M)^2,
\end{equation}
and
\begin{equation}\label{local estimate-4}
  \int_0^{t_0}\big\|\big(\theta_{kx}^{(n)},\chi_{kx}^{(n)}\big)(\tau)\big\|_2^2d\tau\leq C(m,M)(2M)^2.
\end{equation}
Combining with the energy inequality \eqref{local estimate-1}, by direct computation on \eqref{iterative format 1}, one derives
\begin{equation}\label{local estimate-5}
 \|v_k^{(n)}(t)-\bar v\|_2^2\leq (2M)^2,\qquad\int_0^{t_0}\|v_{kx}^{(n)}(\tau)\|_1^2d\tau\leq C(m,M)(2M)^2,
\end{equation}
and 
\begin{equation}\label{local estimate-6}
  \inf_{\mathbb{R}\times[0,t_0]}\theta^{(n)}_k(x,t)\geq \frac{m_1}{2}>0,\quad\inf_{\mathbb{R}\times[0,t_0]}v^{(n)}_k(x,t)\geq \frac{m_2+b+h}{2}.
\end{equation}
Using \eqref{local estimate-1}-\eqref{local estimate-6},  one finally obtains 
\begin{equation}\label{local estimate}
  \big(v_k^{(n)},u^{(n)}_k,\theta^{(n)}_k,\chi_k^{(n)}\big)\in X_{\frac{m_1}{2},\frac{m_2+b+h}{2},2M}(0,t_0).
\end{equation}

Therefore, we are ready to take the limit of the above approximate solution sequence $\{(v^{(n)}_k,u^{(n)}_k,\theta^{(n)}_k,\chi^{(n)}_k)\}_{n=1}^{+\infty}$ obtained by the iterative method. From \eqref{First iteration}, noticing that for fixed $k$, $\{(v^{(n)}_k,u^{(n)}_k,\theta^{(n)}_k)\}_{n=1}^{+\infty}$, $\{\chi_{kx}^{(n)}\}_{n=1}^{+\infty}$ and $\{(\chi_k^{(n)})^2-1\}_{n=1}^{+\infty}$ can be shown as the Cauchy sequences in $C([0,t_0];H^3(\mathbb{R}))$, $C([0,t_0];H^2(\mathbb{R}))$ and $C([0,t_0];H^1(\mathbb{R}))$ respectively. Thus, a solution  $(v_k,u_k,\theta_k,\chi_k)$  can be obtained by letting $n\rightarrow +\infty$. By the same way, taking $T^*$ smaller than $t_0$ if necessary,  $\{v_k\}_{k=1}^{+\infty}$, $\{(u_k,\theta_k,\chi_{kx})\}_{k=1}^{+\infty}$ and  $\{\chi^2_k-1\}_{k=1}^{+\infty}$ are also the Cauchy sequences in $C([0,T^*];H^2(\mathbb{R}))$, $C([0,T^*];H^1(\mathbb{R}))$ and $C([0,T^*];L^2(\mathbb{R}))$ respectively. Therefore, the desired local solution  $(v,u,\theta,\chi)\in X_{\frac{m_1}{2},\frac{m_2+b+h}{2},2M}(0,T^*)$ for \eqref{NSFAC-Lagrange} is finally obtained by letting $k\rightarrow\infty$. The uniqueness of the local solution can be obtained by the basic energy estimation. The proof of Proposition 2.1 is completed.
\end{proof}

\section{A Priori Estimates for Theorem 1.1.}\label{sec3}
\ \ \ \ In this section, we will extend the local strong solution of \eqref{NSFAC-Lagrange}-\eqref{initial condition}  to the whole interval  $[0,T]$ by providing a priori estimates in  Proposition \ref{a priori estimate}, and then thus obtain Theorem \ref{thm-global}. The completion of this extension can be directly guaranteed by the establishment of Proposition \ref{a priori estimate}. Under the framework of the proof for Proposition \ref{a priori estimate}, the restriction of $\beta>0$  is necessary, because of the nature and requirements of \eqref{eta}. On the other hand, since the Navier-Stokes system \eqref{NSF-Lagrange} can be regarded as a special case of system \eqref{NSFAC-Lagrange}, Proposition \ref{a priori estimate} is also valid for the Navier-Stokes system \eqref{NSF-Lagrange} as  $\beta>0$.

\begin{proposition}\textbf{(A priori estimates)}  \label{a priori estimate}
Under the assumptions of Theorem 1.1, suppose  the Cauchy problem \eqref{NSFAC-Lagrange}-\eqref{initial condition} has a solution $(v,u,\theta,\chi)\in X_{m_1,m_2,M}(0,T)$ for some   $M>0,m_1>0,m_2>b+h$ and $T>0$. Then it holds
\begin{eqnarray}\label{a priori estimate inequality}\left.\begin{array}{llll}
\displaystyle\sup_{0\leq t\leq T}\Big(\big\|\big(v-\bar v,u,\theta-\bar\theta,\chi_x\big)\big\|_1^2+\big\|\chi^2-1\big\|^{2}\Big)\\
\displaystyle\qquad+\int_0^T\Big(\big\|v_x(\tau)\big\|^2+\big\|\big(u_x,\theta_x,\chi_{x}\big)(\tau)\big\|_1^2+\big\|\big(\frac{\chi_{x}}{v}\big)_{xx}(\tau)\big\|^2\Big)d\tau\leq C,
\end{array}\right.
\end{eqnarray}
where the positive constant $C>0$ depends on the initial data $v_0,\chi_0,\theta_0,u_0$,  $T$ and $\epsilon,h$. Moreover,   there exists a positive constant $M_1$, such that 
\begin{equation}\label{Upper and lower bounds on v}
 b+h< v \leq M_1<+\infty,\qquad0<\theta<M_1<+\infty,\qquad |\chi|\leq 1,\qquad\forall(x,t)\in\mathbb{R}\times[0,T].
\end{equation}
\end{proposition}

\

The proof of the Proposition 3.1 can be obtained from the following series of lemmas, which we describe as follows. First, let us define the following two key functions
\begin{equation}\label{Phi}
\Phi(v)\xlongequal{\mathrm{def}}\left\{\begin{array}{llll}
\displaystyle\frac1{\bar\theta}\Big(\big(\frac{R\bar\theta}{\bar v-b}-\frac{a}{\bar v^2}\big) (v-\bar v)-(\frac{a}v-\frac{a}{\bar v})-R\bar\theta\ln\frac{v-b}{\bar v-b}\Big), &\displaystyle\forall v>b+h,\\
\displaystyle+\infty,&\displaystyle \forall v\leq b+h.
 \end{array}\right.
\end{equation}
and
\begin{equation}\label{Psi}
\Psi(\theta)\xlongequal{\mathrm{def}} \frac1{\bar\theta} (\theta-\bar\theta)-\ln \frac{\theta}{\bar\theta},  \qquad \forall\theta>0.
 \end{equation}

\

With the help of the positive properties of  $\Psi$ and $\Phi$ above,    the basic energy estimate is established in the following Lemma \ref{Fundamental energy inequality}. 

\begin{lemma}\label{Fundamental energy inequality} Under the assumptions of Proposition \ref{a priori estimate},   it holds that
\begin{equation}\label{basic energy inequality}
\sup_{0\leq t\leq T}\int_{-\infty}^{+\infty}\Big(\frac{1}{2\bar\theta}u^2+\frac1{\bar\theta}W(\chi,\chi_x)+\Phi(v)+\Psi(\theta)\Big)dx+\int_0^TV(t)dt\leq E_{0},
\end{equation}
where
\begin{equation}\label{V(t)}
  W(\chi,\chi_x)\overset{\mathrm{def}}{=}\frac{1}{4\epsilon}(\chi^2-1)^2+\frac\epsilon2\frac{\chi_x^2}{v},\quad
V(t)\overset{\mathrm{def}}{=}\int_{-\infty}^{+\infty}\Big(\frac{\theta^\beta\theta_x^2}{v\theta^2}+\frac{u_x^2}{v\theta}+\frac{v\mu^2}{\theta}\Big)dx,
\end{equation}
and
\begin{equation}\label{E0}
E_0\overset{\mathrm{def}}{=}\int_{-\infty}^{+\infty}\Big[\frac{1}{2\bar\theta}u_0^2+\frac1{\bar\theta}W(\chi_0,\chi_{0x})+
 \Phi(v_0)+\Psi(\theta_0)\Big]dx.
\end{equation}
Moreover, $\forall n=0,\pm1,\pm2,\cdots$,
there exist points $a_{n}(t)$ and $b_n(t)$ on the interval $(n,n+1)$, such that
\begin{equation}\label{bar v theta}
\begin{array}{llll}
 \displaystyle v(a_n(t),t)\overset{\mathrm{def}}{=}\int_n^{n+1}v(x,t)dx\in[ \alpha_1,\alpha_2], \quad \theta(b_n(t),t)\overset{\mathrm{def}}{=}\int_n^{n+1}\theta(x,t)dx\in[ \alpha_1,\alpha_2], 
\end{array} \end{equation}
 where $0<\alpha_1<\alpha_2$    are respectively the two roots of the following algebraic equation
 \begin{equation}\label{algebraic equation}
  \frac1{\bar\theta}\big(\frac{R}{\bar v-b}-\frac{a}{\bar v^2}\big)(y-\bar v)-R\ln\frac{y-b}{\bar v-b}+ \frac1{\bar\theta}(y-\bar\theta)-\ln \frac{y}{\bar\theta} =E_0+\frac{a}{\bar\theta\bar v(b+h)}\int_{-\infty}^{+\infty}(v_0-\bar v)dx.
 \end{equation}
\end{lemma}
\begin{proof}  From \eqref{NSFAC-Lagrange},  \eqref{initial condition}, one has
\begin{equation}\label{conservation}
 \int_{-\infty}^{+\infty}(v(x,t)-\bar v)dx=\int_{-\infty}^{+\infty}(v_0(x)-\bar v)dx.
\end{equation}
Firstly, by direct calculation,  it can be seen that both of $\Phi(v)$ and $\Psi(\theta)$ are positive. That is to say,
\begin{equation}\label{the sign of Phi}
   \Phi(v)> 0,\qquad \Psi(\theta)>0\qquad \forall v\not=\bar{v},\quad\theta\not=\bar{\theta}.
\end{equation} 
This is because,  the function $\Phi$ can be expressed in the following form of pressure integration
\begin{equation}\label{The integral form of phi}
\displaystyle\Phi(v)=\left\{\begin{array}{llll}
\displaystyle \frac{1}{\bar\theta}\Big(p(\bar v,\bar\theta)(v-\bar v)-\int_{\bar v}^vp(\xi,\bar\theta)d\xi\Big), &\displaystyle\forall v>b+h,\\
\displaystyle+\infty,&\displaystyle \forall v\leq b+h,\end{array}\right.
\end{equation}
and the functions $\Phi(v)$ and $\Psi(\theta)$ have the following basic properties
 \begin{equation}\label{properities of phi and psi}
 \left.\begin{array}{llll}
  \displaystyle \Phi(\bar v)=0,\ \ \Psi(\bar\theta)=0;\\
\displaystyle \Phi'(v)\big|_{v=\bar v}=\frac 1{\bar\theta}\big(\frac{R\bar\theta}{\bar v-b}+\frac{a}{v^2}-\frac{R\bar\theta}{v-b}\big)\big|_{v=\bar v}=0,\quad  \Psi'(\theta)\big|_{\theta=\bar\theta}=\big(\frac1{\bar\theta}-\frac1{\theta}\big)\big|_{\theta=\bar\theta}=0;\\     
\displaystyle \Phi''(v)=-\frac 1{\bar\theta}\big(\frac{2a}{v^3}-\frac{R\bar\theta}{(v-b)^2}\big)>0,\  \forall v\in B(\bar v,\delta);\quad \Psi''(\theta)=\frac{1}{\theta^2}>0,\ \forall\theta\neq 0;
\end{array}
\right.
\end{equation}
where $B(\bar v,\delta)$ represents the $\delta$-neighborhood of $\bar v$.
Therefore, combining with \eqref{bar v, bar theta-1}--\eqref{bar v, bar theta-2}, we know that at this point,  $\bar v$ is in the stable region (see \eqref{ABC} and figure 1.(b))  under the fixed temperature $\bar\theta$, combining with \eqref{The integral form of phi},  then, \eqref{the sign of Phi} is immediately obtained.

Secondly, multiplying \eqref{NSFAC-Lagrange}$_1$ by $\frac1{\bar\theta}\big(\frac{R\bar\theta}{\bar v-b}-\frac{a}{\bar v^2}\big)-\frac{R}{v-b}$, \eqref{NSFAC-Lagrange}$_2$ by $\frac1{\bar\theta}u$, \eqref{NSFAC-Lagrange}$_3$ by $\frac1{\bar\theta}-\frac1{\theta}$, and \eqref{NSFAC-Lagrange}$_4$ by $\frac1{\bar\theta}\mu$,  adding them together, one has
\begin{equation}\label{basic equality}
\left.\begin{array}{llll}
\displaystyle\Big(\frac{1}{2\bar\theta}u^2+\frac{1}{4\bar\theta\epsilon}(\chi^2-1)^2+\frac{\epsilon}{2\bar\theta}\frac{\chi_x^2}{v}+\Phi(v)+\Psi(\theta)\Big)_t+
\frac{\theta^\beta\theta_x^2}{v\theta^2}+\frac{u_x^2}{v\theta}+\frac{v\mu^2}{\theta}\\
\displaystyle=\frac1{\bar\theta}\big(\frac{R\bar\theta}{\bar v-b}-\frac{a}{\bar v^2}\big)
u_x+\frac1{\bar\theta}\Big(\frac{uu_x}{v}-\frac{Ru\theta}{v-b}+\frac{au}{v^2}\Big)_x\\
\displaystyle\qquad\qquad+\Big((\frac1{\bar\theta}-\frac1\theta)\frac{\theta^\beta\theta_x}{v}\Big)_x+\frac{\epsilon}{\bar\theta}\big(\frac{\chi_x\chi_t}{v}\big)_x
-\frac{\epsilon}{2\bar\theta}\big(\frac{\chi_x^2u}{v^2}\big)_x.\end{array}\right.
\end{equation}
Integrating \eqref{basic equality} over $(-\infty,+\infty)\times[0,T]$ by parts, then \eqref{basic energy inequality} is obtained.

Finally, considering the positivity of the second-order derivative as follows
$$  \Big(\frac1{\bar\theta}\big(\frac{R\bar\theta}{\bar v-b}-\frac{a}{\bar v^2}\big)(y-\bar v)-R\ln\frac{y-b}{\bar v-b}+\frac1{\bar\theta}(y-\bar\theta)-\ln \frac{y}{\bar\theta}\Big)''>0,$$ 
utilizing the properties of convex functions and Jensen's  inequality, in combination with the basic energy estimate \eqref{basic energy inequality}, one has
\begin{equation}\label{Jensen}
\alpha_1\leq\int_{n}^{n+1}v dx\leq\alpha_2,\quad \alpha_1\leq\int_{n}^{n+1}\theta dx\leq\alpha_2,
\end{equation}
 where
 $0<\alpha_1<\alpha_2$ are two roots of \eqref{algebraic equation}, which implies 
  \begin{equation}\label{roots}
\alpha_1\leq v(a_n(t),t)\leq\alpha_2,\quad \alpha_1\leq\theta(a_n(t),t)\leq\alpha_2.   
 \end{equation}
The proof of Lemma \ref{Fundamental energy inequality} is finished.
\end{proof}

\

In order to obtain the upper bound of $v$, with the inspiration from  Kazhikhov \cite{K1982} and the help of Lemma \ref{Fundamental energy inequality}, the expression of $v$ \eqref{expression of v}  will be established in the following Lemma \ref{lem-expression of v}.

\begin{lemma}\label{lem-expression of v}
Under the assumptions of Proposition \ref{a priori estimate},    then $\forall x\in[n,n+1), n=0,\pm1,\pm2,\cdots$, it has the following expression of $v$
\begin{equation}\label{expression of v}
  v(x,t)=   D(x,t) Y(t)\Big(1+\int_0^t\frac{\big(-\frac{a}{v}+\frac{Rv\theta}{v-b}+\frac{\epsilon}{2}\frac{\chi^2_x}{v}\big)(x,\tau)}{ D(x,\tau) Y(\tau)}d\tau\Big),
\end{equation}
 where
 \begin{equation}\label{D}
  D(x,t)=v_0(x)e^{\int_{n}^{x}(u(y,t)-u_0(y))dy},
   \end{equation}
 and
 \begin{equation}\label{Y}
  Y(t)=\frac{v(n,t)}{v_0(n)}e^{-\int_0^t(-\frac{a}{v^2}+\frac{R\theta}{v-b}+\frac{\epsilon}{2}\frac{\chi^2_x}{v^{2}})(n,s)ds}.
 \end{equation}
\end{lemma}
\begin{proof} Substituting \eqref{NSFAC-Lagrange}$_1$ into \eqref{NSFAC-Lagrange}$_2$, one has
\begin{align}
\displaystyle (\ln v)_{xt} =\big(-\frac{a}{v^2}+\frac{R\theta}{v-b}+\frac{\epsilon}{2}\frac{\chi_x^2}{v^2}\big)_x+u_t,
\end{align}
integrating the above equation over $(0,t)$, one derives
\begin{align}\label{lnv1}
\displaystyle (\ln v)_{x} =\big(\int_0^t(-\frac{a}{v^2}+\frac{R\theta}{v-b}+\frac{\epsilon}{2}\frac{\chi_x^2}{v^2})d\tau\big)_x+u-u_0+(\ln v_0)_{x}.
\end{align}
By using the definition of $a_n(t)$  in \eqref{bar v theta}, for $\forall x\in[n,n+1),n=0,\pm1,\pm2,\cdots$, integrating \eqref{lnv1} from $n$ to $x$ by parts, one gets
\begin{equation}\label{v}
 v(x,t) = D(x,t) Y(t) e^{\int_0^t\big(-\frac{a}{v^2}+\frac{R\theta}{v-b}+\frac{\epsilon}{2}\frac{\chi^2_x}{v^2}\big)(x,s)ds},
 \end{equation}
where $D(x,t)$ and $Y(t)$ are as defined in  \eqref{D} and \eqref{Y}  respectively.  Since
\begin{equation}\label{bound of u}
  \Big|\int_n^x\big(u(y,t)-u_0(y)\big)dy\Big|\leq\Big(\int_n^{n+1}u^2dy\Big)^\frac12+\Big(\int_n^{n+1}u_0^2dy\Big)^\frac12\leq C,
\end{equation}
one immediately obtains
\begin{equation}\label{upper and lower bound of D}
  C^{-1}\leq D(x,t)=v_0(x)e^{\int_{n}^{x}(u(y,t)-u_0(y))dy}\leq C, \ \forall x\in[n,n+1),n=0,\pm1,\pm2,\cdots,
\end{equation}
where $C$ is a constant independent of $n,x$. Now we introduce the function 
\begin{equation}\label{g} 
 g(x,t) \overset{\mathrm{def}}{=}\int_0^t\big(-\frac{a}{v^2}+\frac{R\theta}{v-b}+\frac{\epsilon}{2}\frac{\chi^2_x}{v^2}\big)(x,s)ds,
\end{equation}
combining with \eqref{v}, one  gets the following ordinary differential equation for $g(x,t)$
\begin{equation*} 
 g_t=\frac{-\frac{a}{v}+\frac{Rv\theta}{v-b}+\frac{\epsilon}{2}\frac{\chi^2_x(x,t)}{v(x,t)}}{v(x,t)}
=\frac{-\frac{a}{v}+\frac{R\theta v}{v-b}+\frac{\epsilon}{2}\frac{\chi^2_x(x,t)}{v(x,t)}}{ D(x,t) Y(t)e^g},
\end{equation*}
which gives
\begin{equation*}
  e^g  =1+\int_0^t\frac{-\frac{a}{v}+\frac{Rv\theta}{v-b}+\frac{\epsilon}{2}\frac{\chi^2_x(x,\tau)}{v(x,\tau)}}{ D(x,\tau) Y(\tau)}d\tau,
  \end{equation*}
substituting  the expression above into \eqref{v},  one obtains \eqref{expression of v}.    Thus the proof of Lemma \ref{lem-expression of v} is finished.
\end{proof}
 
 \
 
In the following Lemma \ref{the lower bounds of density and temperature},  the upper and lower bounds of the phase field $\chi$ will be presented through the energy estimation method. Furthermore, by using the truncation technique introduced by Jiang \cite{J1999} and Li-Liang \cite{LL2016} (see\eqref{bar v theta} and \eqref{the truncation function method}),  the lower bound of temperature $\theta$ will be given.
\begin{lemma}\label{the lower bounds of density and temperature}
Under the assumptions of Proposition \ref{a priori estimate},    it holds that 
 \begin{equation}\label{lower bound of density and theta}
v(x,t)> b+h,\quad|\chi(x,t)|\leq C, \quad \theta(x,t)\geq C, \quad\forall(x,t)\in(-\infty,+\infty)\times [0,T],\\
  \end{equation}
where the positive constant $C$ only depends on $T$, the initial data $v_0,u_0,\theta_0,\chi_0$ and $\epsilon,h$.
\end{lemma}
\begin{proof}
First of all, we shall prove that $v$ has a positive lower bound, and this can be directly derived from the definition of van der Waals equation of state \eqref{modified van der Waals equation} and the basic energy estimation. In fact, from \eqref{Phi} and \eqref{basic energy inequality},  one obtains 
  \begin{equation}\label{the lower bound of v}
   v(x,t)> b+h,\qquad \forall(x,t)\in(-\infty,+\infty)\times[0,T],
  \end{equation}
and the direct calculation yields
\begin{equation}\label{v-inequality}
\frac{1}{v-b}<\frac{1}{h}.
\end{equation}

Secondly,  we will consider providing the upper and lower bounds of $\chi$. For $\forall n=0,\pm1,\pm2,\cdots$,  by using \eqref{basic energy inequality}, one has
  \begin{eqnarray}
   \int_n^{n+1}\chi^4(x,t)dx &=& \int_n^{n+1}(\chi^2(x,t)-1+1)^2dx\notag\\\\
   &\leq& 2\int_n^{n+1}(\chi^2(x,t)-1)^2dx+2\notag\\
   &\leq& C,\notag
 \end{eqnarray}
 which yields
  \begin{equation}\label{phi L4}
   \int_n^{n+1}\chi^4(x,t)dx\leq C,
 \end{equation}
then one gets
 \begin{equation}\label{phi L1}
   \int_n^{n+1}\chi(x,t) dx\leq C,
 \end{equation}
where $C$ is independent of $n$. Therefore, combining with \eqref{basic energy inequality}, $\forall (x,t)\in[n,n+1)\times[0,T]$, one obtains
\begin{eqnarray}\label{lower bound of phi}
|\chi(x,t)|&\leq&\big|\int_n^{n+1}\chi(x,t)-\chi(y,t)dx\big|+\big|\int_n^{n+1}\chi(y,t)dy\big|\notag \\
 &\leq&\big|\int_n^{n+1}\big(\int_y^x\chi_z(z,t)dz\big)dy\big| +C\notag\\
 &\leq&\displaystyle\big|\int_n^{n+1}\big(\int_y^x\frac{\chi_z(z,t)}{v^\frac12}v^\frac12dz\big)dy\big| +C\\
&\leq& \big(\int_n^{n+1}\frac{\chi_x^2}{v}dx\big)^{\frac12}\big(\int_{n}^{n+1}vdx\big)^{\frac12}+C\notag \\
 &\leq& \displaystyle C_1,\notag
\end{eqnarray}
where $C_1$ is independent of $n$, thus, the upper and lower bounds of $\chi$ is obtained. 

Finally,  we will provide the lower bound of $\theta$. Denoting by
\begin{equation}\label{set of theta}
(\theta>2\bar\theta)(t)\overset{\mathrm{def}}{=}\big\{x\in\mathbb{R}\big|\theta(t)>2\bar\theta\big\},\quad (\theta<\frac12\bar\theta)(t)\overset{\mathrm{def}}{=}\big\{x\in\mathbb{R}\big|\theta(t)<\frac12\bar\theta\big\}.
\end{equation}
By using \eqref{basic energy inequality}, one has
\begin{eqnarray*}
  E_0&\geq&\int_{(\theta<\frac12\bar\theta)(t)}\Psi(\theta)dx+\int_{(\theta>2\bar\theta)(t)}\Psi(\theta)dx \\
  &\geq&\big(\ln2-\frac12\big)\big|(\theta<\frac12\bar\theta)\big| +\big(1-\ln2\big)\big|(\theta>2\bar\theta)\big|, 
\end{eqnarray*}
which leading to the following inequality\begin{equation}\label{upper bound of the set theta} 
 \big|(\theta<\frac12\bar\theta)\big| +\big|(\theta>2\bar\theta)\big|\leq C.
 \end{equation}
Denoting
\begin{equation}\label{the truncation function method}
  \big(\theta^{-1}-2\bar\theta\big)_+\overset{\mathrm{def}}{=}\max\big\{\theta^{-1}-2\bar\theta,0\big\},
\end{equation}
for $\forall p>2$, multiplying \eqref{NSFAC-Lagrange}$_3$ by $\theta^{-2}\big(\theta^{-1}-2\bar\theta\big)_+^p$, integrating over $(-\infty,+\infty)$  with respect to $x$, combining with \eqref{the lower bound of v} and \eqref{upper bound of the set theta},  one has
\begin{align}\label{lower bound of theta}
&\frac{1}{p+1}\frac{d}{dt}\int_{-\infty}^{+\infty}\big( {\theta}^{-1}-2\bar\theta\big)_+^{p+1} dx+\int_{-\infty}^{+\infty}\big(\frac{u_x^2}{v\theta^2}+\frac{v\mu^2}{\theta^2}\big) \big({\theta}^{-1}-2\bar\theta\big)_+^{p}dx\notag\\
&\qquad+2\int_{-\infty}^{+\infty}\frac{\theta^\beta\theta_x^2}{v\theta^3}\big( {\theta}^{-1}-2\bar\theta\big)_+^{p}dx+p\int_{-\infty}^{+\infty}\frac{\theta^\beta\theta_x^2}{v\theta^4}\big( {\theta}^{-1}-2\bar\theta\big)_+^{p-1}dx\notag\\
 &\leq\int_{-\infty}^{+\infty}\frac{R(1+\frac bh)u_x}{v\theta}\big( {\theta}^{-1}-2\bar\theta\big)_+^{p} dx\\ 
 & \leq\frac{1}{2}\int_{-\infty}^{+\infty}\frac{u_x^2}{v\theta^2}\big( {\theta}^{-1}-2\bar\theta\big)_+^{p}dx+C\int_{-\infty}^{+\infty}\frac{1}{v}\big( {\theta}^{-1}-2\bar\theta\big)_+^{p}dx\notag\\ 
 &\leq\frac{1}{2}\int_{-\infty}^{+\infty}\frac{u_x^2}{v\theta^2}\big( {\theta}^{-1}-2\bar\theta\big)_+^{p} dx+C\Big(\int_{-\infty}^{+\infty}\big( {\theta}^{-1}-2\bar\theta\big)_+^{p+1}\Big)^{\frac{p}{p+1}}.\notag\end{align}
Applying Gronwall's inequality to  the above result \eqref{lower bound of theta}, one derives
\begin{equation}\label{lower bound of theta Lp norm}  \sup_{0\leq t\leq T}\left\|  \big( {\theta}^{-1}-2\bar\theta\big)_+(\cdot,t) \right\|_{L^{p+1}(\mathbb{R})}\leq C,\ \ \forall p>2,
\end{equation}
where $C$ is independent of $p$, and further, taking the limit as $p\rightarrow+\infty$ on both sides of  \eqref{lower bound of theta Lp norm}, one  eventually gets 
\begin{equation}\label{lower bound of theta L-infty norm} 
 \sup_{0\leq t\leq T}\left\|  \big( {\theta}^{-1}-2\bar\theta\big)_+(\cdot,t) \right\|_{L^{\infty}(\mathbb{R})}\leq C,
\end{equation}
which the positive lower bound of temperature $\theta$ is thus proved, and the proof of Lemma 3.3 is completed.
\end{proof}

\

Now,  based on the aforementioned several lemmas, the upper bound estimate of $v$ can be obtained as follows. It is worth noting that the selection of $\eta$ in \eqref{eta} is crucial, which enables us to obtain better regularity regarding the temperature $\theta$, thus together with the help of the key inequality \eqref{The square modulus of the first derivative for phi}, allows us to obtain the another key inequality \eqref{integral sup theta}, and ultimately arrive at the upper bound for $v$ as follows.

\begin{lemma}\label{the upper bound of density}
Under the assumptions of Proposition \ref{a priori estimate},   $\forall T>0$,  it holds that 
 \begin{equation}\label{upper bound of density}
v(x,t)\leq C, \quad\forall(x,t)\in(-\infty,+\infty)\times [0,T],\\
  \end{equation}
where the positive constant $C$ only depends on $T$, the initial data $v_0,u_0,\theta_0,\chi_0$ and $\epsilon,h$.
\end{lemma}
\begin{proof} First, deriving from \eqref{NSFAC-Lagrange}$_5$,  one has
\begin{equation}\label{second derivative form for phi}
 \epsilon\Big(\frac{\chi_x}{v}\Big)_x=-\mu+\frac{1}{\epsilon}(\chi^3-\chi),
\end{equation}
and then combining  \eqref{basic energy inequality}, \eqref{V(t)}, \eqref{conservation}, Lemma \ref{the lower bounds of density and temperature} and \eqref{lower bound of phi} yields
\begin{equation}\label{estimate for second derivative form of phi}
\int_{-\infty}^{+\infty}\Big(\frac{\chi_x}{v}\Big)_x^2\frac{1}{\theta}dx\leq C\big(1+V(t)\big),
\end{equation}
thus, combining with \eqref{basic energy inequality},\eqref{conservation},\eqref{the lower bound of v},\eqref{lower bound of phi}, \eqref{estimate for second derivative form of phi}, one derives
\begin{align}\label{The square modulus of the first derivative for phi}
 \sup_{x\in\mathbb{R}}\big(\frac{\chi_x}{v}\big)^2(x,t)& \leq C\int_{-\infty}^{+\infty}\frac{\chi_x}{v}\Big(\frac{\chi_x}{v}\Big)_xdx\notag\\
 &\leq C\int_{-\infty}^{+\infty}\frac{\theta\chi_x^2}{v^2}dx+\int_{-\infty}^{+\infty}\Big(\frac{\chi_x}{v}\Big)_x^2\frac{1}{\theta}dx\\
  &\leq C\big(\sup_{x\in\mathbb{R}} \theta+1+V(t)\big).\notag
\end{align}

Next, setting
\begin{equation}\label{eta}
  \eta\overset{\mathrm{def}}{=}\frac{1}{4}\max\{1,2-\beta\}\in\big(\frac14,\frac{1}{2}\big),
  \end{equation}
 from Cauchy inequality, combining with \eqref{basic energy inequality}, one has
\begin{align}\label{theta0}
&\int_0^T\sup_{x\in\mathbb{R}}\theta dt\leq C\int_0^T\sup_{x\in\mathbb{R}}\int_{-\infty}^{x}\partial_y\big(\theta-(2\bar\theta)\big)_+dydt+C\notag \\
&\leq C\int_0^T\int_{(\theta>2\bar\theta)(t)}\frac{\theta^\beta \theta_x^2}{v\theta^{2-2\eta}}dxdt+C\int_0^T\int_{(\theta>2\bar\theta)(t)}\frac{v\theta^{2-2\eta}}{\theta^\beta} dxdt +C\\
&\leq C\int_0^T\int_{(\theta>2\bar\theta)(t)}\frac{\theta^\beta \theta_x^2}{v\theta^{2-2\eta}}dxdt+C\int_0^T\sup_{x\in\mathbb{R}}\theta^{\max\{2-2\eta-\beta,0\}}\int_{(\theta>2\bar\theta)(t)}v dxdt +C\notag\\
&\leq C\int_0^T\int_{(\theta>2\bar\theta)(t)}\frac{\theta^\beta \theta_x^2}{v\theta^{2-2\eta}}dxdt+\delta\int_0^T\sup_{\mathbb{R}}\theta dt+C(\delta),\notag
\end{align}
and it immediately follows that
\begin{equation}\label{The integral of theta}
  \int_0^T\sup_{x\in\mathbb{R}}\theta dt\leq \frac{C}{1-\delta}\int_0^T\int_{(\theta>2\bar\theta)(t)}\frac{\theta^\beta \theta_x^2}{v\theta^{2-2\eta}}dxdt+C(\delta).
\end{equation}
Multiplying \eqref{NSFAC-Lagrange}$_3$ by $\big(\theta^{\eta}-(2\bar\theta)^{\eta}\big)_+\theta^{\eta-1}$ and integrating over $(-\infty,+\infty)\times[0,T]$ by parts, combining  \eqref{upper bound of the set theta} and \eqref{basic energy inequality}, then 
\begin{align*}
 & (1-2\eta)\int_0^T\int_{(\theta>2\bar\theta)(t)}\frac{\theta^\beta\theta_x^2}{v\theta^{2-2\eta}}dxdt+\int_0^T\int_{-\infty}^{+\infty}\Big(\frac{u_x^2}{v}+v\mu^2\Big)
\Big(\theta^{\eta}-(2\bar\theta)^{\eta}\Big)_+\theta^{\eta-1}dx\nonumber \\
 &=\frac{1}{2\eta}\int_{-\infty}^{+\infty}\Big(\big(\theta^{\eta}-(2\bar\theta)^{\eta}\big)_+^2-\big(\theta_0^{\eta}-(2\bar\theta)^{\eta}\big)_+^2\Big)dx
 +(2\bar\theta)^\eta(1-\eta)
 \int_0^T\int_{(\theta>2\bar\theta)(t)}\frac{\theta^\beta\theta_x^2}{v\theta^{2-\eta}}dxdt\nonumber\\
 &\qquad+\int_0^T\int_{-\infty}^{+\infty}\frac{R\theta u_x}{v-b}\Big(\theta^{\eta}-(2\bar\theta)^{\eta}\Big)_+\theta^{\eta-1}dxdt\\
&\leq C+\frac{1-2\eta}2\int_0^T\int_{(\theta>2\bar\theta)(t)}\frac{\theta^\beta\theta_x^2}{v\theta^{2-2\eta}}dxdt+
\frac12\int_0^T\int_{-\infty}^{+\infty}\frac{u_x^2}{v}\Big(\theta^{\eta}-(2\bar\theta)^{\eta}\Big)_+\theta^{\eta-1}dxdt  \nonumber\\
&\qquad+C\int_0^T\int_{-\infty}^{+\infty}\frac{\theta^2}{v}\Big(\theta^{\eta}-(2\bar\theta)^{\eta}\Big)_+\theta^{\eta-1}dx dt,\nonumber
\end{align*}
which yields 
\begin{equation}\label{auxiliary inequality for estimation of temperature 1}
\left.\begin{array}{llll}
\displaystyle\frac{1-2\eta}{2}\int_0^T\int_{(\theta>2\bar\theta)(t)}\frac{\theta^\beta\theta_x^2}{v\theta^{2-2\eta}}dxdt+\int_0^T\int_{-\infty}^{+\infty}\Big(\frac{u_x^2}{2v}+v\mu^2\Big)
\displaystyle\Big(\theta^{\eta}-(2\bar\theta)^{\eta}\Big)_+\theta^{\eta-1}dx \\
\displaystyle\leq C\int_0^T\int_{(\theta>2\bar\theta)(t)}\theta^{\eta+1}\Big(\theta^{\eta}-(2\bar\theta)^{\eta}\Big)dx dt+C \\
\displaystyle\leq C\int_0^T\sup_{\mathbb{R}}\theta^{2\eta}\int_{(\theta>2\bar\theta)(t)}\theta dx dt+C\int_0^T\sup_{\mathbb{R}}\theta^{2\eta}\int_{(\theta>2\bar\theta)(t)}\theta^{1-\eta} dx dt+C\\
\displaystyle \leq C\int_0^T\sup_{x\in\mathbb{R}}\theta^{2\eta} dt+C\leq\delta \int_0^T\sup_{x\in\mathbb{R}}\theta dt+C(\delta).
\end{array}\right.
\end{equation}
Thus, combining \eqref{auxiliary inequality for estimation of temperature 1} and \eqref{The integral of theta}, one  derives
\begin{equation}\label{integral sup theta}
  \int_0^T\sup_{x\in\mathbb{R}}\theta dt+\int_0^T\int_{(\theta>2\bar\theta)(t)}\frac{\theta^\beta\theta_x^2}{v\theta^{2-2\eta}}dxdt\leq C.
\end{equation}
and moreover, together with \eqref{basic energy inequality}, one obtains
\begin{equation}\label{integral sup theta 2}
  \int_0^T\int_{-\infty}^{+\infty}\frac{\theta^\beta\theta_x^2}{v\theta^{2-2\eta}}dxdt\leq C.
\end{equation}

Finally, since \eqref{expression of v} yields
\begin{equation}
    Y^{-1}(t)\int_n^{n+1}v(x,t)dx= \int_n^{n+1}  D(x,t)\Big(1+ \int_0^t \frac{\big(-\frac{a}{v}+\frac{Rv\theta}{v-b}+\frac{\epsilon}{2}\frac{\chi^2_x}{v}\big)(x,\tau)}{ D(x,\tau) Y(\tau)}d\tau\Big)dx,\notag
\end{equation}
 together with \eqref{bar v theta}, \eqref{Y}, \eqref{the lower bound of v}, \eqref{upper and lower bound of D} and\eqref{v-inequality}, and combining with Gronwall inequality, there exists  a positive constant $C$, satisfying
\begin{equation}\label{Y1}
 C^{-1}\leq Y^{-1}(t)\leq C+C\int_0^tY^{-1}(s)ds,
 \end{equation}
which implies that
\begin{equation}\label{upper and lower bound for Y}
  0<C^{-1}\leq Y(t)\leq C<+\infty,\qquad \forall(x,t)\in(-\infty,+\infty)\times [0,T].
\end{equation}
Thus, applying  \eqref{upper and lower bound of D}, \eqref{the lower bound of v}, \eqref{v-inequality}, \eqref{upper and lower bound for Y} to \eqref{expression of v}, one obtains
\begin{eqnarray}\label{the upper bound estimate for v 1}
  v(x,t)&=& D(x,t) Y(t)\Big(1+ \int_0^t\frac{\big(-\frac{a}{v}+\frac{Rv\theta}{v-b}+\frac{\epsilon}{2}\frac{\chi^2_x}{v}\big)(x,\tau)}{ D(x,\tau) Y(\tau)}d\tau\Big)\notag\\
&\le& C+C\int_0^t\Big(\sup_{x\in\mathbb{R}} \theta (x,\tau)+\sup_{x\in\mathbb{R}}\big(\frac{\chi_x(x,\tau)}{v(x,\tau)}\big)^2\sup_{x\in\mathbb{R}}v(x,\tau)\Big)d\tau.
\end{eqnarray}
 Substituting \eqref{The square modulus of the first derivative for phi} and \eqref{integral sup theta} into \eqref{the upper bound estimate for v 1}, one achieves
 \begin{eqnarray}\label{the upper bound estimate for v}
  v(x,t)\leq C+C\int_0^t\big(\sup_{x\in\mathbb{R}} \theta+1+V(t)\big)\sup_{x\in\mathbb{R}} v(x,\tau)d\tau.
\end{eqnarray}
Therefore, by applying the Gronwall inequality to  \eqref{the upper bound estimate for v}, one eventually obtains
\begin{equation}\label{the upper bound of v}
 v(x,t)\leq C,\qquad\forall (x,t)\in (-\infty,+\infty)\times [0,T],
\end{equation}
and the proof of Lemma \ref{the upper bound of density} is completed.
\end{proof}

\

After obtaining the upper and lower bounds of $v$ as well as the lower bound of temperature
$\theta$, we now begin to present here the high-order estimates of $u, v, \theta$ and $\chi$, in order to prepare
for eventually obtaining the upper bound of temperature $\theta$. Firstly, the  boundedness  of $\|u_x\|_{L^2(0,T;L^2(\mathbb{R}))}$ and $\|\frac{\theta_x}{\sqrt{\theta}}\|_{L^2(0,T;L^2(\mathbb{R}))}$ will be presented in Lemma 3.5 as following.

\begin{lemma}
Under the assumptions of Proposition \ref{a priori estimate},   $\forall T>0$,  it holds that 
 \begin{equation}\label{ux2andthetax2}
\int_0^T\int_{-\infty}^{+\infty}\big(u_x^2+\frac{\theta_x^2}\theta \big)dxdt\leq C,
 \end{equation}
 where the positive constant $C$ only depends on $T$, the initial data $v_0,u_0,\theta_0,\chi_0$ and $\epsilon,h$.
 \end{lemma}
 \begin{proof}
 Multiplying \eqref{NSFAC-Lagrange}$_2$ by $u$ and integrating the resultant with respect of $x$ over $\mathbb{R}$, by using \eqref{basic energy inequality}, \eqref{bound of u},\eqref{lower bound of density and theta}, \eqref{The square modulus of the first derivative for phi} and \eqref{integral sup theta}, one has
 \begin{equation}\label{thetax2}
 \left.\begin{array}{llll}
 \displaystyle \frac12\frac{d}{dt}\int_{-\infty}^{+\infty}u^2dx+\int_{-\infty}^{+\infty}\frac{u_x^2}vdx\\
 \displaystyle=\int_{-\infty}^{+\infty}\big(-\frac a{v^2}+\frac{R\theta}{v-b}\big)u_xdx+\frac\epsilon2\int_{-\infty}^{+\infty}\frac{\chi_x^2}{v^2}u_xdx\\
\displaystyle \leq C\int_{-\infty}^{+\infty}\big(\frac 1{v}+\frac{\theta}{v}\big)|u_x|dx+\int_{-\infty}^{+\infty}\frac{\chi_x^4}{v^3}dx+\frac14\int_{-\infty}^{+\infty}\frac{u_x^2}vdx\\
\displaystyle \leq C\int_{-\infty}^{+\infty}\frac{|\theta-\bar\theta|}{v}|u_x|dx+C\int_{-\infty}^{+\infty}\frac{|v-\bar v|}{v}|u_x|dx +\int_{-\infty}^{+\infty}\frac{\chi_x^4}{v^3}dx+\frac14\int_{-\infty}^{+\infty}\frac{u_x^2}vdx\\ 
\displaystyle \leq C\int_{-\infty}^{+\infty}(\theta-\bar\theta)^2dx+C\int_{-\infty}^{+\infty}(v-\bar v)^2dx+\int_{-\infty}^{+\infty}\frac{\chi_x^4}{v^3}dx+\frac12\int_{-\infty}^{+\infty}\frac{u_x^2}vdx\\
\displaystyle\leq C+C\int_{(\theta>2\bar\theta)(t)}\theta^2dx+C\big(\sup_{x\in\mathbb{R}} \theta+1+V(t)\big)+\frac12\int_{-\infty}^{+\infty}\frac{u_x^2}vdx\\
\displaystyle\leq C+C\big(\sup_{x\in\mathbb{R}} \theta+V(t)\big)+\frac12\int_{-\infty}^{+\infty}\frac{u_x^2}vdx.
 \end{array}
 \right.
 \end{equation}
 Integrating \eqref{thetax2} over $[0,T]$, combining with \eqref{basic energy inequality} and \eqref{integral sup theta}, one gets
 
 \begin{equation}\label{ux2}
   \int_0^T\int_{-\infty}^{+\infty}u_x^2dxdt\leq C,
 \end{equation}
 further, combining with \eqref{integral sup theta 2} and \eqref{lower bound of density and theta}, one obtains
  \begin{equation}\label{thetax2}
   \int_0^T\int_{-\infty}^{+\infty}\frac{\theta_x^2}{\theta}dxdt\leq 
   c\int_0^T\int_{-\infty}^{+\infty}\frac{\theta^\beta\theta_x^2}{v\theta^{2-2\eta}}\theta^{1-2\eta-\beta}dxdt\leq C.
 \end{equation}
 \end{proof}
 
\

Based on Lemma 3.5, the boundedness  of  $\|v_x\|_{L^2(0,T;L^2(\mathbb{R}))}$, $\|\chi_{xx}\|_{L^2(0,T;L^2(\mathbb{R}))}$ and $\|\chi_{t}\|_{L^2(0,T;L^2(\mathbb{R}))}$ can be provided  as following.

\begin{lemma}\label{the square modulus estimate for density}
Under the assumptions of Proposition \ref{a priori estimate},   $\forall T>0$,  it holds that 
 \begin{equation}\label{bound of the square modulus of the first derivative}
\sup_{0\leq t\leq T}\int_{-\infty}^{+\infty}v_x^2dx\leq C,\qquad \int_0^T\int_{-\infty}^{+\infty}(\chi_{xx}^2+\chi_t^2)dxdt\leq C,
 \end{equation}
 where the positive constant $C$ only depends on $T$, the initial data $v_0,u_0,\theta_0,\chi_0$ and $\epsilon,h$.
 \end{lemma}
\begin{proof} First, substituting \eqref{NSFAC-Lagrange}$_1$ into  \eqref{NSFAC-Lagrange}$_2$, one has
\begin{equation}\label{Another form of the momentum equation}
 \Big(u-\frac{v_x}{v}\Big)_t+\big(\frac{R\theta}{v-b}\big)_x=\Big(\frac{a}{v^2}-\frac\epsilon2\big(\frac{\chi_x}{v}\big)^2\Big)_x.
\end{equation}
Multiplying \eqref{Another form of the momentum equation} by $u-\frac{v_x}{v}$, integrating by parts over $(-\infty,+\infty)\times [0,t]$, one has
\begin{equation}\label{the square modulus estimate for density-1}
\left.\begin{array}{llll}
\displaystyle\frac12\int_{-\infty}^{+\infty}\Big(u-\frac{v_x}{v}\Big)^2(x,\tau)dx\Big|_0^t+\int_0^t\int_{-\infty}^{+\infty}\frac{R\theta v_x^2}{v(v-b)^2}dxd\tau \\
\displaystyle=\int_0^t\int_{-\infty}^{+\infty}\Big[\Big(-\frac{2av_x}{v^3}-\frac{R\theta_x}{v-b}-\epsilon\frac{\chi_x}{v}
\big(\frac{\chi_x}{v}\big)_x\Big)\Big(u-\frac{v_x}{v}\Big)+\frac{R\theta v_xu}{(v-b)^2}\Big]dxd\tau.
\end{array}\right.
\end{equation}
On the one hand, from \eqref{lower bound of density and theta}, one gets
\begin{equation}\label{vx-1}
\left.\begin{array}{llll}
 \displaystyle \Big|\int_0^t\int_{-\infty}^{+\infty}\frac{2av_x}{v^3}\big(u-\frac{v_x}{v}\big)dxd\tau\Big|\\
 \displaystyle\leq C\int_0^t\int_{-\infty}^{+\infty}v_x^2dxd\tau+C\int_0^t\int_{-\infty}^{+\infty}\big(u-\frac{v_x}{v}\big)^2dxd\tau,
\end{array}
\right.
\end{equation}
and by using \eqref{ux2andthetax2}, one also derives
\begin{equation}\label{vx-2}
\left.\begin{array}{llll}
 \displaystyle\Big|\int_0^T\int_{-\infty}^{+\infty}\frac{R\theta_x}{v-b}\big(u-\frac{v_x}{v}\big)dxd\tau\Big|\\
 \displaystyle\leq C\int_0^T\int_{-\infty}^{+\infty}\frac{\theta_x^2}{\theta}dxd\tau+C \int_0^T\int_{-\infty}^{+\infty}\frac{\theta}{v^2}\Big(u-\frac{v_x}{v}\Big)^2dxd\tau\\
\displaystyle\leq C+C\int_0^T\sup_{x\in\mathbb{R}}\theta\int_{-\infty}^{+\infty}\Big(u-\frac{v_x}{v}\Big)^2dxd\tau.
\end{array}
\right.
\end{equation}
On the other hand, by using \eqref{basic energy inequality}, \eqref{lower bound of density and theta}, \eqref{upper bound of density}--\eqref{The square modulus of the first derivative for phi} and \eqref{the upper bound of v}, the following inequalities can be derived
\begin{equation}\label{u-vx/v}
\left.\begin{array}{llll}
 \displaystyle\Big|\int_0^t\int_{-\infty}^{+\infty}\epsilon\frac{\chi_x}{v}\big(\frac{\chi_x}{v}\big)_x\Big(u-\frac{v_x}{v}\Big)dxd\tau\Big|\\
 \displaystyle\leq C\int_0^t\int_{-\infty}^{+\infty}\Big(\big|(\frac{\chi_x}{v})_x\big|^2+\big|\frac{\chi_x}{v}\big|^2\big(u-\frac{v_x}{v}\big)^2\Big)dxd\tau+C\\
 \displaystyle\leq C+C\int_0^t\sup_{x\in\mathbb{R}}\big|\frac{\chi_x}{v}\big|^2\int_{-\infty}^{+\infty}\big(u-\frac{v_x}{v}\big)^2dxd\tau
 +C\int_0^t\int_{-\infty}^{+\infty}(\chi_{x}^{2}v_{x}^{2}+\chi_{xx}^{2})dxd\tau,
\end{array}
\right.
\end{equation}
and
\begin{equation}\label{tv}
\left.\begin{array}{llll}
\displaystyle \int_0^t\int_{-\infty}^{+\infty} \frac{R\theta v_xu}{(v-b)^2}dxd\tau\\
\displaystyle \leq  \frac12\int_0^t\int_{-\infty}^{+\infty}\frac{R\theta v_x^2}{v(v-b)^2}dxd\tau+\frac12\int_0^t\int_{-\infty}^{+\infty}\frac{Ru^2\theta v}{(v-b)^2}dxd\tau\\
\displaystyle\leq \frac12\int_0^t\int_{-\infty}^{+\infty}\frac{R\theta v_x^2}{v(v-b)^2}dxd\tau+C\int_0^t\sup_{\mathbb{R}}\theta d\tau\\
\displaystyle\leq \frac12\int_0^t\int_{-\infty}^{+\infty}\frac{R\theta v_x^2}{v(v-b)^2}dxd\tau+C.
 \end{array}
 \right.
\end{equation}
substituting \eqref{vx-1}-\eqref{tv} into \eqref{the square modulus estimate for density-1}, combining Lemma 3.1--Lemma 3.4, by using Gronwall's inequality, one obtains
\begin{equation}\label{estimate for v-x^2}
\left.\begin{array}{llll}
\displaystyle\sup_{0\leq t\leq T} \int_{-\infty}^{+\infty}\Big(u-\frac{v_x}{v}\Big)^2dx+\int_0^T\int_{-\infty}^{+\infty}\frac{\theta v_x^2}{v(v-b)^2}dxdt\\
\displaystyle\leq C+C\int_0^T\int_{-\infty}^{+\infty}(\chi_{x}^{2}v_{x}^{2}+\chi_{xx}^{2})dxdt,
\end{array}\right.
\end{equation}
which together with \eqref{basic energy inequality} implies
 \begin{equation}\label{L-infty v-x}
   \sup_{0\leq t\leq T}\int_{-\infty}^{+\infty}v_x^2dx\leq C+C\int_0^T\int_{-\infty}^{+\infty}(\chi_{x}^{2}v_{x}^{2}+\chi_{xx}^{2})dxdt.
 \end{equation}
Secondly, rewriting  \eqref{NSFAC-Lagrange}$_{3,4}$ as follows
\begin{equation}\label{Allen-Cahn-Lagrange}
 \chi_t-\epsilon\chi_{xx}=-\epsilon\frac{\chi_xv_x}{v}-\frac{v}{\epsilon}\big(\chi^3-\chi\big),
\end{equation}
multiplying \eqref{Allen-Cahn-Lagrange} by $\chi_{xx}$, integrating the resultant over $(-\infty,+\infty)$ with respect of $x$,  by using \eqref{The square modulus of the first derivative for phi} and \eqref{lower bound of phi}, combining with Lemma 3.1-3.4,  one obtains
\begin{align}
&\quad\frac{1}{2}\frac{d}{dt}\int_{-\infty}^{+\infty}\chi_x^2dx+\epsilon\int_{-\infty}^{+\infty}\chi_{xx}^2dx\notag\\
&=\epsilon\int_{-\infty}^{+\infty}\frac{\chi_xv_x}{v}\chi_{xx}dx+ \frac{1}{\epsilon}\int_{-\infty}^{+\infty}v(\chi^3-\chi)\chi_{xx}dx\notag\\
& \leq C(\epsilon)\Big(\int_{-\infty}^{+\infty}\chi_x^2v_x^2dx+\int_{-\infty}^{+\infty}\chi^2(\chi^2-1)^2dx\Big)+\frac{\epsilon}{2}\int_{-\infty}^{+\infty}\chi_{xx}^2dx\\
&\leq   C\int_{-\infty}^{+\infty}\chi_{x}^{2}v_{x}^{2}dx+\frac{\epsilon}{2}\int_{-\infty}^{+\infty}\chi_{xx}^2dx+C,\notag
\end{align}
thus, by using \eqref{basic energy inequality} and \eqref{integral sup theta},  one achieves
\begin{equation}\label{phi-xx-L2}
  \int_0^T\int_{-\infty}^{+\infty}\chi_{xx}^2dxdt\leq C+ C\int_0^T\int_{-\infty}^{+\infty}\chi_{x}^{2}v_{x}^{2}dx.
\end{equation}
Multiplying \eqref{L-infty v-x} by $\frac{1}{2C}$ and adding \eqref{phi-xx-L2}, one has
\begin{align}
&\int_{-\infty}^{+\infty}v_{x}^{2}dx+ \int_0^T\int_{-\infty}^{+\infty}\chi_{xx}^2dxdt\notag\\
&\leq C+ C\int_0^T\int_{-\infty}^{+\infty}\chi_{x}^{2}v_{x}^{2}dx\\
&\leq C+ C\int_0^T\sup_{x\in\mathbb{R}}\frac{\chi_{x}^{2}}{v^{2}} \int_{-\infty}^{+\infty}v_{x}^{2}dx,\notag
\end{align}
together with \eqref{The square modulus of the first derivative for phi} and applying Gronwall inequality, one obtains
\begin{align}
\int_{-\infty}^{+\infty}v_{x}^{2}dx+ \int_0^T\int_{-\infty}^{+\infty}\chi_{xx}^2dxdt\leq C.
\end{align}
Finally,  integrating  \eqref{Allen-Cahn-Lagrange}  over $(-\infty,+\infty)$ with respect of $x$, by using \eqref{L-infty v-x},  one obtains
\begin{align}\label{estimate of phi-xx-l2}
\int_{-\infty}^{+\infty}\chi_t^2dx&\leq C\Big(\int_{-\infty}^{+\infty}\chi_{xx}^2dx+\int_{-\infty}^{+\infty}\chi_x^2v_x^2dx+\int_{-\infty}^{+\infty}\big(\chi^3-\chi\big)^2dx\Big)\notag\\
 &\leq C\Big(\int_{-\infty}^{+\infty}\chi_{xx}^2dx+\int_{-\infty}^{+\infty}v_x^2dx\int_{-\infty}^{+\infty}\chi_{xx}^2dx+1\Big)\\
 &\leq C\Big(\int_{-\infty}^{+\infty}\chi_{xx}^2dx+1\Big),\notag
\end{align}
then, from \eqref{phi-xx-L2}, we achieve
\begin{equation}\label{phi-t-L2}
  \int_0^T\int_{-\infty}^{+\infty}\chi_t^2dxdt\leq C.
\end{equation}
The proof of Lemma \ref{the square modulus estimate for density} is completed.
\end{proof}

\

Due to the appearance of the free energy density $(\frac{\chi_x^2}{v^2})_x$ in the momentum equation, the higher regularity of $\chi$ is necessary, and the specific estimation is as follows.
\begin{lemma}\label{phi-xxx}
Under the assumptions of Proposition \ref{a priori estimate},   $\forall T>0$,  it holds that
 \begin{equation}\label{bound of the square modulus of the third derivative for phi}
\sup_{0\leq t\leq T}\int_{-\infty}^{+\infty}\chi_{xx}^2dx+\int_0^T\int_{-\infty}^{+\infty}\Big(\chi_{xt}^2+\big(\frac{\chi_{x}}{v}\big)_{xx}^2\Big)dxdt\leq C,
 \end{equation}
 where the positive constant $C$ only depends on $T$, the initial data $v_0,u_0,\theta_0,\chi_0$ and $\epsilon,h$.
\end{lemma}
\begin{proof}
Rewriting  \eqref{Allen-Cahn-Lagrange} as follows
\begin{equation}\label{Allen-Cahn-Lagrange-1}
\frac{\chi_t}{v}-\epsilon\Big(\frac{\chi_x}{v}\Big)_{x}=-\frac{1}{\epsilon}(\chi^3-\chi),
\end{equation}
differentiating \eqref{Allen-Cahn-Lagrange-1} with respect to $x$, one has
\begin{equation}\label{Allen-Cahn-Lagrange-2}
 \Big(\frac{\chi_x}{v}\Big)_t-\epsilon\Big(\frac{\chi_x}{v}\Big)_{xx}=-\frac{1}{\epsilon}\big(\chi^3-\chi\big)_x+\frac{\chi_tv_x}{v^2}-\frac{\chi_x u_x}{v^2},
\end{equation}
multiplying \eqref{Allen-Cahn-Lagrange-2} by $\big(\frac{\chi_{x}}{v}\big)_t$, integrating the resultant over $(-\infty,+\infty)$ with respect of $x$, by using \eqref{basic energy inequality}, \eqref{lower bound of density and theta}, \eqref{upper bound of density}, \eqref{integral sup theta},  \eqref{estimate of phi-xx-l2}, one obtains
\begin{eqnarray}
&&\int_{-\infty}^{+\infty}\Big(\frac{\chi_{x}}{v}\Big)_t^2dx+ \frac{\epsilon}{2}\frac{d}{dt}\int_{-\infty}^{+\infty}\Big(\frac{\chi_x}{v}\Big)_{x}^2dx \notag\\
&& =-\frac{1}{\epsilon}\int_{-\infty}^{+\infty}\big(\chi^3-\chi\big)_x\big(\frac{\chi_{x}}{v}\big)_tdx+\int_{-\infty}^{+\infty}\frac{\chi_tv_x}{v^2}\big(\frac{\chi_{x}}{v}\big)_tdx-
\int_{-\infty}^{+\infty}\frac{\chi_x u_x}{v^2}\big(\frac{\chi_{x}}{v}\big)_tdx\notag\\
&&\leq C\Big(\int_{-\infty}^{+\infty}\big(3\chi^2-1\big)^2\chi_x^2dx+\int_{-\infty}^{+\infty}\chi_t^2v_x^2dx+\int_{-\infty}^{+\infty}\chi_x^2 u_x^2dx\Big)+\frac{1}{3}\int_{-\infty}^{+\infty}\big(\frac{\chi_{x}}{v}\big)_t^2dx \notag\\
&&\leq C\Big(1+\|\chi_t\|_{L^\infty}^2\int_{-\infty}^{+\infty}v_x^2dx+\|\frac{\chi_{x}}{v}\|_{L^\infty}^2\int_{-\infty}^{+\infty}u_x^2dx\Big)
+\frac{1}{3}\int_{-\infty}^{+\infty}\big(\frac{\chi_{x}}{v}\big)_t^2dx\\
&&\leq C\Big(1+\int_{-\infty}^{+\infty}\big(\chi_t^2+2|\chi_t\chi_{xt}|\big)dx+\int_{-\infty}^{+\infty}\Big(\frac{\chi_x}{v}\Big)_{x}^2dx\int_{-\infty}^{+\infty}u_x^2dx\Big)
+\frac{1}{3}\int_{-\infty}^{+\infty}\big(\frac{\chi_{x}}{v}\big)_t^2dx\notag\\
&&\leq C\Big(1+\int_{-\infty}^{+\infty}\chi_{t}^2dx+\varepsilon\int_{-\infty}^{+\infty}\chi_{xt}^2dx+\int_{-\infty}^{+\infty}\Big(\frac{\chi_x}{v}\Big)_{x}^2dx\int_{-\infty}^{+\infty}
u_x^2dx\Big)
+\frac{1}{3}\int_{-\infty}^{+\infty}\big(\frac{\chi_{x}}{v}\big)_t^2dx\notag\\
&&\leq C\Big(1+\int_{-\infty}^{+\infty}\chi_{t}^2dx+\int_{-\infty}^{+\infty}u_x^2dx\int_{-\infty}^{+\infty}\Big(\frac{\chi_x}{v}\Big)_{x}^2dx\Big)
+\frac{1}{2}\int_{-\infty}^{+\infty}\big(\frac{\chi_{x}}{v}\big)_t^2dx,\notag
\end{eqnarray}
where $\chi_{xt}=\big(\frac{\chi_x}{v}\big)_t v+\frac{\chi_x u_x}{v}$ is used in the last inequality. Thus, combining with Gronwall's inequality, one derives
\begin{equation}\label{higher derivative energy estimation for phi}
 \sup_{t\in[0,T]} \int_{-\infty}^{+\infty}\Big(\frac{\chi_x}{v}\Big)_{x}^2dx+\int_0^T\int_{-\infty}^{+\infty}\Big(\frac{\chi_{x}}{v}\Big)_t^2dxdt\leq C,
\end{equation}
together with  \eqref{basic energy inequality}, \eqref{bound of the square modulus of the first derivative},  one achieves
\begin{equation}
\sup_{0\leq t\leq T}\int_{-\infty}^{+\infty}\chi_{xx}^2dx+\int_0^T\int_{-\infty}^{+\infty}\chi_{xt}^2dxdt\leq C.
 \end{equation}
Furthermore,   combining with \eqref{basic energy inequality} and \eqref{higher derivative energy estimation for phi}, one obtains
\begin{equation}\label{the upper and lower bounds of the derivative for phi}
  \sup_{(x,t)\in\mathbb{R}\times[0,T]}\Big|\frac{\chi_x}{v}\Big|\leq \sqrt{2}\Big\|\frac{\chi_x}{v}\Big\|\Big\|\big(\frac{\chi_x}{v}\big)_{x}\Big\|\leq C.
\end{equation}
Further, by using \eqref{Allen-Cahn-Lagrange-2} and the energy estimation obtained above, combining with Lemma 3.1--3.5, one achieves
\begin{equation}
\int_0^T\int_{-\infty}^{+\infty}\big(\frac{\chi_{x}}{v}\big)_{xx}^2dxdt\leq C,
 \end{equation}
  the proof of Lemma \ref{phi-xxx} is completed.
\end{proof}

\

After obtaining the  boundedness  of the interface free energy $\|(\frac{\chi_x}{v})_{xx}\|_{L^2(0,T;L^2(\mathbb{R}))}$ ,  the high-order estimates of $u$  will be presented as following.
\begin{lemma}\label{velocity-x}
Under the assumptions of Proposition \ref{a priori estimate},   $\forall T>0$,  it holds that
 \begin{equation}\label{energy estimate of velocity}
\sup_{0\leq t\leq T}\int_{-\infty}^{+\infty}u_{x}^2dx+\int_0^T\int_{-\infty}^{+\infty}\big(u_{t}^2+u^2_{xx}\big)dxdt\leq C,
 \end{equation}
where the positive constant $C$ only depends on $T$, the initial data $v_0,u_0,\theta_0,\chi_0$ and $\epsilon,h$.
\end{lemma}
\begin{proof}
Multiplying \eqref{NSFAC-Lagrange}$_2$ by $u_{xx}$, integrating the resultant over $(-\infty,+\infty)\times(0,T)$ by parts,  combining with \eqref{lower bound of density and theta},  \eqref{bound of the square modulus of the first derivative}, \eqref{integral sup theta}, \eqref{higher derivative energy estimation for phi} and \eqref{the upper and lower bounds of the derivative for phi},  one has
\begin{equation}\label{u-x and u-xx 1}
\left.\begin{array}{llll}
 \displaystyle \frac{1}{2}\int_{-\infty}^{+\infty}u_x^2dx+\int_0^T\int_{-\infty}^{+\infty}\frac{u_{xx}^2}{v}dxdt  \\
 \displaystyle\leq C+\frac{1}{2}\int_0^T\int_{-\infty}^{+\infty}\frac{u_{xx}^2}{v}dxdt\\
 \displaystyle\qquad\qquad+C\int_0^T\int_{-\infty}^{+\infty}\Big(v_x^2+\theta_x^2+\theta^2 v_x^2+\big|\frac{\chi_x}{v}\big|^2\big|\big(\frac{\chi_x}{v}\big)_x\big|^2+u^2_xv_x^2\Big)dxdt\\
\displaystyle\leq C+\frac{1}{2}\int_0^T\int_{-\infty}^{+\infty}\frac{u_{xx}^2}{v}dxdt+C\int_0^T\int_{-\infty}^{+\infty}\theta_x^2dxdt+C\int_0^T\sup_{x\in\mathbb{R}}\theta^2\int_{-\infty}^{+\infty}v_x^2dxdt\\
\displaystyle\ \ +C\sup_{\mathbb{R}\times[0,T]}\big|\frac{\chi_x}{v}\big|^2\int_0^T\int_{-\infty}^{+\infty}\big|\big(\frac{\chi_x}{v}\big)_x\big|^2dxdt+
 C\int_0^T\sup_{x\in\mathbb{R}}u_x^2\int_{-\infty}^{+\infty}v_x^2dxdt\\
\displaystyle\leq C+\frac{3}{4}\int_0^T\int_{-\infty}^{+\infty}\frac{u_{xx}^2}{v}dxdt+C\int_0^T\int_{-\infty}^{+\infty}\frac{\theta^\beta\theta_x^2}{v}dxdt,
\end{array}
\right.
\end{equation}
where the following inequality  are used
\begin{eqnarray}\label{u-x L2}
  \sup_{\mathbb{R}} u_x^2&\leq&\int_{-\infty}^{+\infty}\big|(u_x^2)_x\big|dx\notag\\
  &\leq& 2\Big(\int_{-\infty}^{+\infty}u_x^2dx\Big)^{\frac12}\Big(\int_{-\infty}^{+\infty}u_{xx}^2dx\Big)^{\frac12}\\
  &\leq& C(\delta)\int_{-\infty}^{+\infty}u_x^2dx+\delta\int_{-\infty}^{+\infty}\frac{u_{xx}^2}{v}dxdt,\notag
\end{eqnarray}
and
\begin{equation}\label{theta-x L2}
\left.\begin{array}{llll}
\displaystyle \sup_{\mathbb{R}} (\theta-2\bar\theta)_+^2&=\displaystyle \sup_{\mathbb{R}}\Big(\int_x^{+\infty}\partial_y(\theta-2\bar\theta)_+(y,t)dy\Big)^2\\
\displaystyle  &\displaystyle\leq \Big(\int_{(\theta>2\bar\theta)(t)}\big|\theta_y\big|dy\Big)^2\\
 \displaystyle &\displaystyle\leq C\int_{-\infty}^{+\infty}\theta_x^2dx.
 \end{array}\right.
\end{equation}
Multiplying \eqref{NSFAC-Lagrange}$_3$ by $(\theta-2\bar\theta)_+$, integrating the resultant over $\mathbb{R}\times(0,T)$ by parts, combining with \eqref{basic energy inequality}, \eqref{integral sup theta} and \eqref{u-x L2}, one has
\begin{align}\label{theta theta-x 1}
 & \frac{1}{2}\int_{-\infty}^{+\infty}(\theta-2\bar\theta)_+^2dx+\int_0^T\int_{(\theta>2\bar\theta)(t)}\frac{\theta^\beta\theta_x^2}{v}dxdt \notag \\
 & \leq C\int_0^T\int_{-\infty}^{+\infty}\theta(\theta-2\bar\theta)_+|u_x|dxdt+C\int_0^T\int_{-\infty}^{+\infty}\big(u_x^2+\mu^2\big)(\theta-2\bar\theta)_+ dxdt+C\\
 &\leq C\int_0^T\sup_{\mathbb{R}}\theta\Big(\int_{-\infty}^{+\infty}(\theta-2\bar\theta)_+^2dx+\int_{-\infty}^{+\infty}u_x^2 dx\Big)dt+C.\notag
 \end{align}
Since
\begin{equation}
\left.\begin{array}{llll}
\displaystyle  \int_0^T\int_{-\infty}^{+\infty}\theta^\beta\theta_x^2dxdt&=&\displaystyle\int_0^T\int_{(\theta>2\bar\theta)(t)}\theta^\beta\theta_x^2dxdt
  +\int_0^T\int_{(\theta\leq2\bar\theta)(t)}\theta^\beta\theta_x^2dxdt \\
   &\leq& \displaystyle C \int_0^T\int_{(\theta>2\bar\theta)(t)}\frac{\theta^\beta\theta_x^2}{v}dxdt
  +C\int_0^T\int_{(\theta\leq2\bar\theta)(t)}\frac{\theta^\beta\theta_x^2}{v\theta^2}dxdt \\
   &\leq&\displaystyle C \int_0^T\int_{(\theta>2\bar\theta)(t)}\frac{\theta^\beta\theta_x^2}{v}dxdt+C,
\end{array}
\right.
\end{equation}
together  \eqref{u-x and u-xx 1} and \eqref{theta theta-x 1}, combining with Gronwall's inequality, one obtains
\begin{equation}\label{tempeture and velocity-1}
  \sup_{t\in[0,T]}\int_{-\infty}^{+\infty}\big((\theta-2\bar\theta)_+^2+u_x^2\big)dx+\int_0^T\int_{-\infty}^{+\infty}\big(\theta^\beta\theta_x^2+ u_{xx}^2\big)dxdt\leq C.
\end{equation}
Rewriting \eqref{NSFAC-Lagrange}$_2$ as
\begin{equation}\label{momentum equation}
  u_t=-\Big(-\frac3{v^2}+\frac{\theta}{3v-1}\Big)_x+\frac{u_{xx}}{v}-\frac{u_x v_x}{v^2}-\epsilon\frac{\phi_x}{v}\big(\frac{\phi_x}{v}\big)_x,
\end{equation}
from \eqref{bound of the square modulus of the first derivative}, \eqref{integral sup theta}, \eqref{the upper and lower bounds of the derivative for phi}, \eqref{u-x L2}, \eqref{bound of the square modulus of the third derivative for phi} and \eqref{tempeture and velocity-1}, one gets
\begin{eqnarray}\label{sup u-t L2}
 \int_0^T\int_{-\infty}^{+\infty}u_t^2dxdt&\leq& C\int_0^T\int_{-\infty}^{+\infty}\Big(u_{xx}^2+u_x^2v_x^2++v_x^2+\theta_x^2+\theta^2v_x^2+\big|\frac{\phi_x}{v}\big|^2\big|\big(\frac{\phi_x}{v}\big)_x\big|^2\Big)dxdt\notag\\
 &\leq& C,\notag
\end{eqnarray}
combining with \eqref{tempeture and velocity-1}, the energy inequality \eqref{energy estimate of velocity} is obtained. The proof of Lemma \ref{velocity-x} is completed.
\end{proof}

\

Finally,  the upper bound of the temperature $\theta$ and its higher-order estimates will be given as following, in order to complete the proof of Proposition \ref{a priori estimate}.
\begin{lemma}\label{theta-x}
Under the assumptions of Proposition \ref{a priori estimate},   $\forall T>0$,  it holds that
 \begin{equation}\label{energy estimate of tempeture}
\sup_{(x,t)\in\mathbb{R}\times[0,T]}\theta(x,t)+\sup_{0\leq t\leq T}\int_{-\infty}^{+\infty}\theta_{x}^2dx+\int_0^T\int_{-\infty}^{+\infty}\big(\theta_{t}^2+\theta^2_{xx}\big)dxdt\leq C,
 \end{equation}
where the positive constant $C$ only depends on $T$, the initial data $v_0,u_0,\theta_0,\chi_0$ and $\epsilon,h$.
\end{lemma}
\begin{proof}
Multiplying \eqref{NSFAC-Lagrange}$_3$ by $\theta^\beta\theta_t$,  integrating the resultant over $(-\infty,+\infty)$ with respect of $x$, by using \eqref{lower bound of density and theta}, \eqref{the upper bound of density}, \eqref{tempeture and velocity-1}, one has
\begin{equation}\label{theta-t theta-xx}
\left.\begin{array}{llll}
 \displaystyle\frac{1}{2}\frac{d}{dt}\Big(\int_{-\infty}^{+\infty}\frac{(\theta^\beta\theta_x)^2}{v}dx\Big)+\int_{-\infty}^{+\infty}\theta^\beta\theta_t^2dx \\
\displaystyle =-\frac{1}{2}\int_{-\infty}^{+\infty}\frac{(\theta^\beta\theta_x)^2u_x}{v^2}dx+\int_{-\infty}^{+\infty}\frac{\theta^\beta\theta_t\big(u_x^2+v^2\mu^2\big)}{v}dx-\int_{-\infty}^{+\infty}
  \frac{R\theta^{\beta+1}\theta_t u_x}{v-b}dx\\
 \displaystyle\leq C\sup_{x\in\mathbb{R}}|u_x|\int_{-\infty}^{+\infty}\theta^{2\beta}\theta_x^2dx+\frac{1}{2}\int_{-\infty}^{+\infty}\theta^\beta\theta_t^2dx
\\
\displaystyle\qquad\qquad +C\int_{-\infty}^{+\infty}\theta^{\beta+2}u_x^2dx+C\int_{-\infty}^{+\infty}\theta^\beta\big(u_x^4+\mu^4\big)dx\\
 \displaystyle\leq C\big(1+\|u_{xx}\|^2)\int_{-\infty}^{+\infty}\theta^{2\beta}\theta_x^2dx+\frac{1}{2}\int_{-\infty}^{+\infty}\theta^{\beta}\theta_t^2dx
  +C\sup_{x\in\mathbb{R}}\Big(\theta^{2\beta+2}+u_x^4+\mu^4\Big)+C.
  \end{array}\right.
 \end{equation}
Now we deal with the term $\sup_{x\in\mathbb{R}}\big(\theta^{2\beta+2}+u_x^4+\mu^4\big)$ in the last inequality of \eqref{theta-t theta-xx}.
Applying Lemma \ref{velocity-x}, direct computation shows that
\begin{eqnarray}
  \int_0^T\sup_{x\in\mathbb{R}}u_x^4dt  &\leq& C\int_0^T\int_{-\infty}^{+\infty}|u_x^3u_{xx}|dxdt\notag\\
  &\leq&C\int_0^T\sup_{x\in\mathbb{R}}u_x^2\Big(\int_{-\infty}^{+\infty}u_x^2dx\Big)^\frac{1}{2}\Big(\int_{-\infty}^{+\infty}u_{xx}^2dx\Big)^\frac{1}{2}dt\notag\\
   &\leq&\frac12\int_0^T\sup_{x\in\mathbb{R}}u_x^4dt+C\int_0^T\int_{-\infty}^{+\infty}\big(u_x^2+u_{xx}^2)dxdt\notag\\
      &\leq&\frac12\int_0^T\sup_{x\in\mathbb{R}}u_x^4dt+C,\notag
\end{eqnarray}
then
\begin{equation}\label{u-x^4}
 \int_0^T\sup_{x\in\mathbb{R}}u_x^4dt\leq C.
\end{equation}
By the same way, combining with \eqref{bound of the square modulus of the third derivative for phi},  one gets
\begin{equation}\label{mu^4}
 \int_0^T\sup_{x\in\mathbb{R}}\mu^4dt\leq C.
\end{equation}
Moreover, by employing the Sobolev embedding theorem alongside \eqref{tempeture and velocity-1}, one derives
\begin{align}\label{the upper bound estimate for theta}
  \sup_{x\in\mathbb{R}}\theta^{2\beta+2}
  &\leq C\sup_{\mathbb{R}}\Big(\int_{x}^{+\infty}\partial_y(\theta-2\bar\theta)_+^{\beta+1}dy\Big)^2+C\notag\\
  &\leq C\int_{(\theta>2\bar\theta)(t)}(\theta^\beta\theta_x)^2dx+C\\
  &\leq C\int_{-\infty}^{+\infty}\frac{(\theta^\beta\theta_x)^2}{v}dx+C.\notag
\end{align}
Substituting \eqref{u-x^4}, \eqref{mu^4}, \eqref{the upper bound estimate for theta} into \eqref{theta-t theta-xx}, combining with Gronwall's inequality, one obtains
\begin{equation}\label{derivative estimation for theta}
\sup_{0\leq t\leq T}\int_{-\infty}^{+\infty}(\theta^\beta\theta_x)^2dx+\int_0^T\int_{-\infty}^{+\infty}\theta^\beta\theta_t^2dxdt\leq C.
\end{equation}
Thus,  together with \eqref{the upper bound estimate for theta}, it yields
\begin{equation}\label{the upper bound estimate for temperature}
 \sup_{(-\infty,+\infty)\times[0,T]}\theta\leq C.
\end{equation}
And therefore, both \eqref{the upper bound estimate for theta} and \eqref{derivative estimation for theta} derive that
\begin{align}
\displaystyle  \sup_{0 \le t\le T}\int_{-\infty}^{+\infty}\theta_{x}^2 dx+\int_0^T\int_{-\infty}^{+\infty} \theta_t^2dxdt\le C.
\end{align}
moreover, it follows from \eqref{basic energy inequality} and \eqref{lower bound of density and theta} that 
\begin{equation}\label{v-bar v, theta-bar theta}
\left.\begin{array}{llll}
\displaystyle\sup_{0\leq t\leq T}\int_{-\infty}^{+\infty}\big((v-\bar v)^2+(\theta-\bar\theta)^2\big)dx\\
\displaystyle\leq C\sup_{0\leq t\leq T}\int_{-\infty}^{+\infty}\Big( \lambda(v-\bar v)-R\ln\frac{v-b}{\bar v-b}+\frac1{\bar\theta}(\theta-\bar\theta)-\ln \frac{\theta}{\bar\theta}\Big)dx\\
\displaystyle\leq C.
\end{array}\right.
\end{equation}
Now, rewriting  \eqref{NSFAC-Lagrange}$_3$ as
\begin{equation}\label{energy conservation equation theta-xx}
\frac{\theta^\beta\theta_{xx}}{v}=\theta_t- \frac{\beta\theta^{\beta-1}\theta_x^2}{v}+\frac{\theta^\beta\theta_x v_x}{v^2}+\frac{R\theta u_x}{v-b}-\frac{u_x^2+(v\mu)^2}{v},
\end{equation}
which yields  that
\begin{align}\label{the estimate of tempture-xx L2}
\int_0^T\int_{-\infty}^{+\infty}\theta_{xx}^2dxdt&\leq C \int_0^T\int_{-\infty}^{+\infty}\big(\theta_x^2v_x^2+\theta_x^4+u_x^4+\mu^4+u_x^2+\theta_t^2\big)dxdt\notag\\
 &\leq C(\delta)+C\delta\int_0^T\sup_{x\in\mathbb{R}}\theta_x^2dt\\
  &\leq C(\delta)+C\delta\int_0^T\int_{-\infty}^{+\infty}\theta_{xx}^2dxdt.\notag
\end{align}
Furthermore, by using maximum principle \cite{p2005}, we obtain $-1\leq\phi\leq1$. 
The proof of Lemma \ref{theta-x} is completed.
\end{proof}

\section{A Priori Estimate for Theorem 1.2.}\label{sec4}
\ \ \ \ In this section, we will present the proof of Theorem \ref{NS}. Since the compressible Navier-Stokes system \eqref{NSF-Lagrange} can be regarded as a special case of system \eqref{NSFAC-Lagrange}, Proposition \ref{local existence} and Proposition \ref{a priori estimate} is also valid for the  system \eqref{NSF-Lagrange}, and the Theorem \ref{NS} can be regarded as a corollary  of Theorem \ref{thm-global} as  $\beta>0$. Therefore,  we only need to consider whether Theorem \ref{NS} holds when $\beta=0$.

In fact, considering that the local solution which presented in  Proposition \ref{local existence} holds for all $\beta\geq0$ , and noting that in the a prior estimation given in Proposition \ref{a priori estimate} , nly the upper bound of $v$ was obtained by using the condition that $\beta > 0$ (see \eqref{eta}), therefore, in reality, we only need to provide an upper bound proof for $v$ when $\beta=0$ for the compressible Navier-Stokes system \eqref{NSF-Lagrange}, that is, we only need to extend Lemma \ref{the upper bound of density} to the case where $\beta = 0$, and the remaining estimates in  Proposition \ref{a priori estimate} are obviously straightforward. For the sake of the article's brevity, we will not elaborate on this here. 

 Regarding the generalization of Lemma \ref{the upper bound of density} when $\beta = 0$,  noting that \eqref{eta} no longer holds true as $\beta = 0$, and whether \eqref{eta} holds or not becomes the key to obtaining the upper bound of $v$ from the expression \eqref{expression of v} according to the proof framework of the section 3. Therefore, we revised the proof strategy for the upper bound of $v$ in Section 3. Specifically, we choose to start with the estimation of $\big\|\frac{v_x}{v}\big\|^2_{L^2(\mathbb{R})}$, and then use the inequality \eqref{M-v} to obtain the upper bound of $v$. For this purpose, in addition to the basic energy estimate which is already obtained in Lemma \ref{Fundamental energy inequality}, through delicate energy estimates, we further discovered some new and more precise energy inequalities \eqref{theta L-1}-\eqref{theta-bartheta} about $v$ and $\theta$.  From this, we can obtain the upper bound estimate for $v$ through the new energy estimation method we have adopted here.  The specific proof is as follows.
\begin{lemma}\label{beta=0}
Under the assumptions of \eqref{modified van der Waals equation}--\eqref{Hypothesis of p}, \eqref{bar v, bar theta-1}-\eqref{bar v, bar theta-2}, the thermal conductivity coefficient $\kappa$ is a positive constant, and
\begin{equation}\label{condition 1'}(v_0-\bar v)\in L^1(\mathbb{R}),\
(v_0-\bar v,u_0)\in H^1(\mathbb{R}),\quad \theta_0-\bar \theta\in H^1(\mathbb{R}),\ \inf_{x\in\mathbb{R}}v_0(x)>0,\  \inf_{x\in\mathbb{R}}\theta_0(x)>0,
\end{equation}
 it holds that
 \begin{equation}\label{upper bound of density-2}
v(x,t)\leq C,
 \end{equation}
where the positive constant $C$ only depends on $T$, the initial data $v_0,u_0,\theta_0,\chi_0$ and $\epsilon,h$.
\end{lemma}
\begin{proof} It is noted that in the proof of the Section 3, except for the inequalities related to \eqref{eta}, no requirement for $\beta > 0$ is imposed. therefore, the Lemma \ref{Fundamental energy inequality}  still holds for \eqref{NSF-Lagrange} when $\beta \geq 0$, that is
\begin{equation}\label{basic energy inequality for NSF}
\sup_{0\leq t\leq T}\int_{-\infty}^{+\infty}\Big(\frac{1}{2\bar\theta}u^2+\Phi(v)+\Psi(\theta)\Big)dx+\int_0^T\int_{-\infty}^{+\infty}\Big(\frac{\theta_x^2}{v\theta^2}+\frac{u_x^2}{v\theta}\Big)dxdt\leq E_{0},
\end{equation}
where $\Phi$ and $\Psi$ are defined by \eqref{Phi} and \eqref{Psi}. Denoting
\begin{equation}\label{tilde V}
  \tilde V(t)=\int_{-\infty}^{+\infty}\frac{u^2}{2\bar\theta}dx+\int_{-\infty}^{+\infty}\Big(\frac{\theta_x^2}{v\theta^2}+\frac{u_x^2}{v\theta}\Big)dx+1.
\end{equation}
Obviously, one has
\begin{equation}\label{upper bound of tilde V}
\int_0^T\tilde V(t)dt\leq C.
\end{equation}
To establish a rigorous upper bound for $\sup_{\mathbb{R}} v$,  we define
\begin{equation}\label{M-v}
 M_v(t)=1+\max_{\mathbb{R}}v(x,t).
\end{equation}
Before presenting the argument, we first establish several key inequalities that will be used in energy estimation in the proof. By using \eqref{upper bound of the set theta}, and \eqref{basic energy inequality},  the following energy estimates can be derived:
\begin{equation}\label{theta L-1}
 \int_{(\theta\leq2\bar\theta)(t)}(\theta-\bar\theta)^2dx+\int_{(\theta>2\bar\theta)(t)}\theta dx\leq\int_{-\infty}^{+\infty}\Psi(\theta)dx\leq C,
\end{equation}
and
\begin{equation}\label{ln v}
 \int_{-\infty}^{+\infty}\ln^2 \frac{v-b}{\bar v-b} dx\leq\int_{(v\leq2\bar v)(t)}(v-\bar v)^2dx+C\int_{(v>2\bar v)(t)}vdx\leq C\int_{-\infty}^{+\infty}\Phi(v)dx\leq C.
\end{equation}
Moreover, by using \eqref{ln v}, noting that
\begin{eqnarray}\label{v-2}
  (v-2\bar v)_+&\leq&C\Big(\int_{(v>2\bar v)(t)}v^2dx\Big)^{\frac12}\Big(\int_{-\infty}^{+\infty}\frac{v_x^2}{v^2}dx\Big)^{\frac12}\notag\\
   &\leq& C \Big(M_v(t)\int_{(v>2)(t)}vdx\Big)^{\frac12}\Big(\int_{-\infty}^{+\infty}\frac{v_x^2}{v^2}dx\Big)^{\frac12}\\
   &\leq&C M_v^{\frac12}(t)\Big(\int_{-\infty}^{+\infty}\frac{v_x^2}{v^2}dx\Big)^{\frac12},\notag
\end{eqnarray}
one derives that
\begin{equation}\label{Mv(t)}
 M_v(t)\leq C+C\int_{-\infty}^{+\infty}\frac{v_x^2}{v^2}dx.
\end{equation}
Moreover, one can also obtains
\begin{equation}\label{theta-bartheta}
\left.\begin{array}{llll}
\displaystyle \int_0^t\int_{-\infty}^{+\infty}(\theta-\bar\theta)^2dxd\tau\\
\displaystyle \leq  \int_0^t\int_{(\theta<2\bar\theta)(t)}(\theta-\bar\theta)^2dxd\tau+C \int_0^t\int_{-\infty}^{+\infty}(\theta^{\frac12}-\bar\theta^{\frac12})_+^2\theta dxd\tau\\
\displaystyle \leq \int_0^t\max_{\mathbb{R}}(\theta^{\frac12}-\bar\theta^{\frac12})_+^2\int_{(\theta>2\bar\theta)(t)}\theta dxd\tau+C\\
\displaystyle\leq \int_0^t\Big(\int_{(\theta>2\bar\theta)(t)}|\theta_x|\theta^{-\frac12} dx\Big)^2d\tau+C\\
 \displaystyle \leq \int_0^t\Big(\int_{-\infty}^{+\infty}\frac{\theta_x^2}{\theta^2 v}dx\Big)\Big(\int_{(\theta>2\bar\theta)(t)}v\theta dx\Big)d\tau+C\\
 \displaystyle \leq\int_0^t \tilde V(\tau)M_v(\tau)d\tau+C.
\end{array}
\right.
\end{equation}

Based on the above preparatory work, we will now present the proof of the Lemma \ref{beta=0}. Observing that
\begin{equation}\label{identical relation-1}
 \big(\frac{u_x}{v}\big)_x=\big(\frac{v_t}{v}\big)_x=\big(\ln v\big)_{xt}=\big(\frac{v_x}{v}\big)_t,
\end{equation}
then \eqref{NSF-Lagrange}$_2$ can be rewritten as
\begin{equation}\label{the other form for NSAH-2}
 \big(\frac{v_{x}}{v}\big)_t- \big(p'_v(v,\theta))v_x+p_\theta'(v,\theta)\theta_x\big)=u_t.
\end{equation}
Multiplying  \eqref{the other form for NSAH-2} by $\frac{v_x}{v}$, integrating over $(-\infty,+\infty)$ by parts, combining with \eqref{identical relation-1}, one has
\begin{equation}\label{the basic energy equality-2 for density}
\left.\begin{array}{llll}
\displaystyle\frac d{dt}\int_{-\infty}^{+\infty}\Big(\frac{1}{2}\big(\frac{v_x}{v}\big)^2-u\frac{v_x}{v}\Big)dx+\int_{-\infty}^{+\infty}\frac{R\theta}{(v-b)^2}\frac{v_x^2}{v}dx\\
\displaystyle
=\int_{-\infty}^{+\infty}\Big(\frac {2a}{v^2}\frac{v_x^2}{v^2}+\frac{R\theta_x}{v-b}\frac{v_x}v+\frac{u_x^2}{v}\Big)dx.
\end{array}\right.
\end{equation}
Further, multiplying \eqref{NSF-Lagrange}$_2$ by $u$, integrating over $(-\infty,+\infty)\times[0,t]$, by using \eqref{ln v} and \eqref{theta-bartheta}, one has
\begin{equation}\label{the basic energy equality-3-1 for density}
\left.\begin{array}{llll}
\displaystyle
\sup_{[0,t]}\int_{-\infty}^{+\infty}\frac{ u^2}{2}dx+\int_0^t\int_{-\infty}^{+\infty}\frac{ u_x^2}{v}dxd\tau\\
\displaystyle\leq\int_0^t\int_{-\infty}^{+\infty}\big(\frac{R\theta}{v-b}-\frac{a}{v^2}\big)u_xdxd\tau+C\\
\displaystyle\leq C\int_0^t\int_{-\infty}^{+\infty}\Big(\Big|\frac{(\theta-\bar\theta)u_x}{v}\Big|+\Big|\frac{(v-\bar v)u_x}{v}\Big|+
\Big|\frac{v^{\frac12}}{(v-\bar v)^{\frac12}}(v-\bar v)^{\frac12}\frac{u_x}{v^{\frac12}}\Big|\Big)dxd\tau+C\\
\displaystyle \leq \frac14\int_0^t\int_{-\infty}^{+\infty}\frac{u_x^2}{v}dxd\tau+C\int_0^t\int_{-\infty}^{+\infty}(\theta-\bar\theta)^2dxd\tau+C\int_0^t\int_{-\infty}^{+\infty}(v-\bar v)^2dxd\tau+C\\
\displaystyle \leq\frac14\int_0^t\int_{-\infty}^{+\infty}\frac{u_x^2}{v}dxd\tau+C\int_0^t \tilde V(\tau)M_v(\tau)d\tau\\
\displaystyle\ \ \ \ \ +C\int_0^t\Big(\int_{(v\leq 2\bar v)(t)}(v-\bar v)^2dx+M_v(s)\int_{(v>2\bar v)(t)}v dx\Big)d\tau+C\\
\displaystyle\leq \frac14\int_0^t\int_{-\infty}^{+\infty}\frac{u_x^2}{v}dxd\tau+C\int_0^t \tilde V(\tau)M_v(\tau)d\tau+C,
\end{array}\right.
\end{equation}
and then, one gets
 \begin{equation}\label{the basic energy equality-3 for density}
\displaystyle
\sup_{[0,t]}\int_{-\infty}^{+\infty}\frac{ u^2}{2}dx+\frac34\int_0^t\int_{-\infty}^{+\infty}\frac{ u_x^2}{v}dxd\tau\leq C\int_0^t \tilde V(\tau)M_v(\tau)d\tau,
\end{equation}
multiplying \eqref{the basic energy equality-2 for density} by $\frac{1}{2}$ and integrating the resultant over $[0,t]$, adding \eqref{the basic energy equality-3 for density}  together,  it is derived that
\begin{equation}\label{the energy equality for v-x-1}
\left.\begin{array}{llll}
\displaystyle\int_{-\infty}^{+\infty}\Big(\frac{1}{2}u^2-\frac{1}{2}u\frac{v_x}{v}+\frac{1}{4}\big(\frac{v_x}{v}\big)^2\Big)(t)dx+\int_0^t\int_{-\infty}^{+\infty}
\Big(\frac{R\theta}{2(v-b)^2}\frac{v_x^2}{v}+\frac{ u_x^2}{4v}\Big)dxd\tau\\
\displaystyle\leq\int_0^t\int_{-\infty}^{+\infty}\frac{R\theta_x}{2(v-b)}\frac{v_x}vdxd\tau+C\int_0^t \tilde V(\tau)M_v(\tau)d\tau+C\int_0^t\int_{-\infty}^{+\infty} \frac{v_x^2}{v^2}dxd\tau+C.
\end{array}\right.
\end{equation}
Observing that \eqref{NSF-Lagrange}$_3$ can be reformulated as
$$\big(\frac{\theta_x}{v}\big)_x=\theta_t+\frac{R\theta}{v-b}u_x-\frac{u_x^2}{v},$$
integrating the first term of the inequality \eqref{the energy equality for v-x-1} over $[0,t]$ with respect of $\tau$, one has
\begin{equation}\label{I}
\left.\begin{array}{llll}
\displaystyle \int_0^t\int_{-\infty}^{+\infty}\frac{R\theta_x}{2(v-b)}\frac{v_x}vdx= \frac R2\int_0^t\int_{-\infty}^{+\infty}\frac{\theta_x}{v}(\ln \frac{v-b}{\bar v-b})_xdxd\tau\\
 \displaystyle= -\frac R2\int_0^t\int_{-\infty}^{+\infty}\Big(\frac{\theta_x}{v}\Big)_x\ln \frac{v-b}{\bar v-b}dxd\tau \\
\displaystyle= -\frac R2\int_0^t\int_{-\infty}^{+\infty}\Big((\theta-\bar\theta)_t+\frac{R\theta}{v-b}u_x-\frac{u_x^2}{v}\Big)\ln \frac{v-b}{\bar v-b} dxd\tau\\
\displaystyle \leq -\frac R2\Big[\int_{-\infty}^{+\infty}\Big((\theta-\bar\theta)\ln \frac{v-b}{\bar v-b}\Big)(t) dx-\int_0^t\int_{-\infty}^{+\infty}\frac{(\theta-\bar\theta)u_x}{v-b}dxd\tau\\
\displaystyle\qquad\qquad+\int_0^t\int_{-\infty}^{+\infty}\frac{\theta u_x}{v-b}\ln \frac {v-b}{\bar v-b}dxd\tau\Big]+C\int_0^t\Big(\int_{-\infty}^{+\infty}\frac{u_x^2}v dx\Big) M_v d\tau+1\Big)\\
\displaystyle =-\frac R2\sum_{k=1}^3I_{k}+C\Big(\int_0^t\big(\int_{-\infty}^{+\infty}\frac{u_x^2}v dx\big) M_v d\tau+1\Big).
\end{array}
\right.
\end{equation}
For the term $I_1$, combining with \eqref{basic energy inequality}, \eqref{theta L-1}, \eqref{ln v} and \eqref{Mv(t)}, one has
\begin{equation}\label{I1}
\left.\begin{array}{llll}
I_{1}&\leq&\displaystyle\int_{-\infty}^{+\infty}|\theta-\bar\theta|\big|\ln \frac{v-b}{\bar v-b}\big| (t)dx\\
&\leq&\displaystyle \int_{(\theta\leq2\bar\theta)(t)}(\theta (t)-\bar\theta)^2dx+ \int_{-\infty}^{+\infty}\ln^2 \frac{v(t)-b}{\bar v-b} dx+\int_{(\theta>2\bar\theta)(t)}\theta dx\ln M_v(t)\\
&\leq& \displaystyle C+C\ln\int_{-\infty}^{+\infty}(\ln v)_x^2 (t)dx\\
&\leq& \displaystyle C+\delta\int_{-\infty}^{+\infty}(\ln v)_x^2 (t)dx.
\end{array}
\right.
\end{equation}
The estimate for $I_2$ follows from \eqref{theta L-1}, \eqref{theta-bartheta}, \eqref{the basic energy equality-3 for density},
\begin{equation}\label{I2}
\left.\begin{array}{llll}
I_{2} &\leq& \displaystyle\int_0^t\int_{-\infty}^{+\infty}v^{-\frac12}\big|(\theta-\bar\theta)\big|\big|u_xv^{-\frac12}\big|dxd\tau \\
 &\leq&\displaystyle C \int_0^t\int_{-\infty}^{+\infty}\frac{u^2_x}{v}dxd\tau+C \int_0^t\int_{-\infty}^{+\infty}(\theta-\bar\theta)^2dxd\tau\\
 &\leq&\displaystyle C\int_0^t \tilde V(\tau)M_v(\tau)d\tau+C.
\end{array}
\right.
\end{equation}
As for  the term  $I_{3}$, an application of inequalities \eqref{ln v}, \eqref{theta-bartheta} and \eqref{the basic energy equality-3 for density} leads to 
\begin{equation*}
\left.\begin{array}{llll}\label{I3}
I_{3}&=&\displaystyle \int_0^t\int_{-\infty}^{+\infty}\frac{\theta u_x}{v}\ln \frac{v-b}{\bar v-b}dxd\tau\\
&\leq& \displaystyle C\Big(\int_0^t\int_{-\infty}^{+\infty}\Big|\frac{(\theta-\bar\theta) u_x}{v}\ln \frac{v-b}{\bar v-b}\Big|dxd\tau+\int_0^t\int_{-\infty}^{+\infty}\Big|\frac{ u_x}{v}\ln \frac{v-b}{\bar v-b}\Big|dxd\tau\Big)\\
&\leq&\displaystyle \int_0^t\int_{-\infty}^{+\infty} \frac{(\theta-\bar\theta)^2}{v} \ln^2 \frac{v-b}{\bar v-b}dxd\tau+C\int_0^t\int_{-\infty}^{+\infty}\frac{u_x^2}{v}dxd\tau+C\int_0^t\int_{-\infty}^{+\infty}\ln^2 \frac{v-b}{\bar v-b}dxd\tau\\
&\leq&\displaystyle C\int_0^t\int_{-\infty}^{+\infty}\frac{u_x^2}{v}dxd\tau+C\int_0^t\int_{-\infty}^{+\infty}(\theta-\bar\theta)^2dxd\tau+C\\
&\leq&\displaystyle C \int_0^t\tilde V(\tau)M_v(\tau) d\tau+C,
\end{array}\right.
\end{equation*}
where the estimation formula $\frac{1}{v}\ln\frac{v-b}{\bar v-b}\leq C$ is used.
 Plugging $I_{1}$, $I_{2}$ and $I_{3}$ into \eqref{I}, one obtains
\begin{equation}\label{theta-x-lnv-x}
\left.\begin{array}{llll}
\displaystyle\int_0^t\int_{-\infty}^{+\infty}\frac{R\theta_x}{2(v-b)}\frac{v_x}vdx\\
\displaystyle\leq C\int_0^t\Big(\tilde V(\tau)M_v(\tau)+1 \Big) d\tau+C+C\delta\int_{-\infty}^{+\infty}(\ln v)_x^2 (t)dx.
\end{array}
\right.
\end{equation}
Together with \eqref{the energy equality for v-x-1} and \eqref{Mv(t)}, it is derived that
\begin{equation}\label{the energy equality for v-x-2}
\left.\begin{array}{llll}
\displaystyle\sup_{[0,t]}\int_{-\infty}^{+\infty}\Big(\frac{1}{2}u^2-\frac{1}{2}u\frac{v_x}{v}+\frac3{16}\big(\frac{v_x}{v}\big)^2\Big)dx+\int_0^t\int_{-\infty}^{+\infty}
\Big(\frac{R\theta}{2(v-b)^2}\frac{v_x^2}{v}+\frac{ u_x^2}{4v}\Big)d\tau\\
\displaystyle \leq C\int_0^t\tilde V(\tau)M_v(\tau)d\tau+C\int_0^t \int_{_\infty}^{+\infty}\frac{v^2_x}{v^2}d\tau+C,
\end{array}
\right.
\end{equation}
and therefore, by combining the Gronwall inequality, combining with \eqref{upper bound of tilde V} and \eqref{Mv(t)}, it can be immediately concluded that
\begin{equation}\label{upper bound of v-x}
M_v(t)\leq  C.
\end{equation}
The proof of Lemma \ref{beta=0} is completed.
\end{proof}

\

At this point,  the upper bound of $v$ is obtained. Then, following the framework of  Section 3, without any difficulty, we can obtain similar a priori estimates for the Navier-Stokes system \eqref{NSF-Lagrange} as $\beta=0$. Thus, we can extend the local solution to the global solution of \eqref{NSF-Lagrange}-\eqref{initial condition of NSF} succussed. In order to avoid repetition, we will omit the detailed process of this proof here.

\section*{Funding} 
Yazhou Chen acknowledges support from National Natural Sciences Foundation of China (NSFC) (No. 12471207),(No. 12371434).\\
Yi Peng acknowledges support from National Natural Sciences Foundation of China (NSFC) (No. 12301266).\\
Xiaoding Shi acknowledges support from National Natural Sciences Foundation of China (NSFC) (No. 12171024), Beijing Natural Science Foundation  (No. 1252004).\\
Xiaoping Wang acknowledges support from National Natural Science Foundation of China (NSFC) (No. 12271461), the key project of NSFC (No. 12131010), the University Development Fund from The Chinese University of Hong Kong, Shenzhen (UDF01002028), and Hetao Shenzhen-Hong Kong Science and Technology Innovation Cooperation Zone Project (No. HZQSWS-KCCYB-2024016).

\end{document}